\documentclass[aap,MSNbibl,nameyear,seceqn,dvips]{arximspdf}
\usepackage{mathrsfs}
\usepackage{accents}
\usepackage{dcolumn}
\usepackage{graphicx}

\doi{10.1214/11-AAP780}
\volume{22}
\issue{2}
\pubyear{2012}
\firstpage{576}
\lastpage{607}

\makeatletter

\newcolumntype{d}[1]{D{.}{.}{#1}}

\newtheorem{theorem}{Theorem}[section]

\newtheorem{lemma}{Lemma}[section]
\newtheorem{corollary}{Corollary}[section]

\newproclaim{definition}{Definition}[section]
\newproclaim{approximation}{Approximation}[section]

\newcommand{\sA}{\mathscr{A}}
\newcommand{\sL}{\mathscr{L}}

\newcommand{\bbE}{\mathbb{E}}
\newcommand{\bbN}{\mathbb{N}}
\newcommand{\bbP}{\mathbb{P}}
\newcommand{\bbR}{\mathbb{R}}

\newcommand{\cP}{{\mathcal P}}

\newcommand{\gd}{\delta}
\newcommand{\gr}{\rho}
\newcommand{\gs}{\sigma}
\newcommand{\bfmath}[1]{\mathbf{#1}}
\newcommand{\bfmathh}[1]{\bolds{#1}}

\newcommand{\bfa}{{\bfmath{a}}}
\newcommand{\bfb}{{\bfmath{b}}}
\newcommand{\bfc}{{\bfmath{c}}}
\newcommand{\bfe}{{\bfmath{e}}}
\newcommand{\bfn}{{\bfmath{n}}}
\newcommand{\bfq}{{\bfmath{q}}}
\newcommand{\bfr}{{\bfmath{r}}}
\newcommand{\bfv}{{\bfmath{v}}}

\newcommand{\bfM}{\bfmath{M}}
\newcommand{\bfP}{\bfmath{P}}
\newcommand{\bfX}{\bfmath{X}}

\newcommand{\bfzero}{\bfmath{0}}
\newcommand{\phammer}[2]{(#1)_{#2}}
\newcommand{\pd}[1]{\frac{\partial}{\partial #1}}
\newcommand{\pdd}[2]{\frac{\partial^2}{\partial #1 \,\partial #2}}
\newcommand{\ud}{\,\mathrm{d}}

\newcommand{\tA}{\theta_{\A}}
\newcommand{\tB}{\theta_{\B}}
\newcommand{\A}{A}
\newcommand{\B}{B}

\newcommand{\td}{\widetilde{d}}
\newcommand{\tbfX}{\widetilde{\bfX}}

\newcommand{\qOTR}[1]{q^{(#1)}_{\mathrm{OTR}}}
\newcommand{\qPS}[1]{q^{(#1)}_{\mathrm{PS}}}
\newcommand{\qaPS}[1]{\qApprox^{(#1)}_{\mathrm{PS}}}
\newcommand{\qPade}[1]{q^{[#1]}_{\mathrm{Pad}\mbox{\fontsize{8.36pt}{8.36pt}\selectfont{\textup{\'e}}}}}
\newcommand{\Pade}{{\mathrm{Pad}\mbox{\fontsize{8.36pt}{8.36pt}\selectfont{\textup{\'e}}}}}
\newcommand{\pade}[3]{[#1]_{#2}(#3)}
\newcommand{\ts}{{\mathrm{s}}}
\renewcommand{\mathring}[1]{\accentset{\circ}{#1}}
\newcommand{\qApprox}{\mathring{q}}
\newcommand{\eApprox}{\mathring{e}}
\newcommand{\ePS}[1]{e_{\mathrm{PS}}^{(#1)}}
\newcommand{\ePade}[1]{e_{\mathrm{Pad}\mbox{\fontsize{8.36pt}{8.36pt}\selectfont{\textup{\'e}}}}^{(#1)}}
\newcommand{\eaPS}[1]{\eApprox_{\mathrm{PS}}^{(#1)}}
\newcommand{\eaPade}[1]{\eApprox_{\mathrm{Pad}\mbox{\fontsize{8.36pt}{8.36pt}\selectfont{\textup{\'e}}}}^{(#1)}}
\newcommand{\abc}{\bfa,\bfb,\bfc}
\newcommand{\floor}[1]{\lfloor #1\rfloor}
\newcommand{\ceil}[1]{\lceil #1\rceil}
\newcommand{\partition}[1]{\cP_{#1}}

\makeatother

\begin{document}
\begin{frontmatter}

\title{Pad\'e approximants and exact two-locus sampling distributions}
\runtitle{Analytic two-locus sampling distributions}

\begin{aug}
\author[A]{\fnms{Paul A.} \snm{Jenkins}\thanksref{t1}\ead[label=e1]{pauljenk@eecs.berkeley.edu}}
and
\author[B]{\fnms{Yun S.} \snm{Song}\corref{}\thanksref{t1,t2}\ead[label=e2]{yss@stat.berkeley.edu}}
\runauthor{P. A. Jenkins and Y. S. Song}
\affiliation{University of California, Berkeley}
\address[A]{Computer Science Division\\
University of California, Berkeley\\
Berkeley, California 94720\\
USA\\
\printead{e1}}
\address[B]{Department of Statistics \\
\quad and Computer Science Division\\
University of California, Berkeley\\
Berkeley, California 94720\\
USA\\
\printead{e2}}
\end{aug}

\thankstext{t1}{Supported in part by NIH Grants R00-GM080099 and R01-GM094402.}
\thankstext{t2}{Supported in part by an Alfred P. Sloan Research
Fellowship and a Packard Fellowship for
Science and Engineering.}

\received{\smonth{7} \syear{2010}}
\revised{\smonth{4} \syear{2010}}

%
\begin{abstract}
For population genetics models with recombination, obtaining an exact,
analytic sampling distribution has remained a challenging open problem
for several decades. Recently, a new perspective based on asymptotic
series has been introduced to make progress on this problem.
Specifically, closed-form expressions have been derived for the first
few terms in an asymptotic expansion of the two-locus sampling
distribution when the recombination rate $\rho$ is moderate to large.
In this paper, a new computational technique is developed for finding
the asymptotic expansion to an \textit{arbitrary order}. Computation in
this new approach can be automated easily. Furthermore, it is proved
here that only a finite number of terms in the asymptotic expansion is
needed to recover (via the method of Pad\'e approximants) the
\textit{exact} two-locus sampling distribution as an analytic function of
$\rho$; this function is exact for all values of $\rho\in[0,\infty)$.
It is also shown that the new computational framework presented here is
flexible enough to incorporate natural selection.
\end{abstract}

%

\pdfkeywords{92D15, 65C50, 92D10, Population genetics,
recombination, sampling distribution,
asymptotic expansion, Pade approximants}

\end{frontmatter}

\section{Introduction}
\label{int}

Central to many applications in genetic analysis is the notion of
sampling distribution, which describes the probability of observing a
sample of DNA sequences randomly drawn from a population. In the
one-locus case with special models of mutation such as the
infinite-alleles model or the finite-alleles parent-independent
mutation model, closed-form sampling distributions have been known
for many decades [\citet{wri1949}, \citet{ewe1972}]. In
contrast, for
multi-locus models with finite recombination rates, finding a
closed-form sampling distribution has so far remained a challenging
open problem. To make progress on this long-standing issue, we
recently proposed a new approach based on asymptotic expansion and
showed that it is possible to obtain useful analytic results when the
recombination rate is moderate to large [Jenkins and Song (\citeyear{jenson2009G},
\citeyear{jenson2010})]. More precisely, our previous work can be
summarized as
follows.

Consider a two-locus Wright--Fisher model, with the two loci denoted by
$A$ and $B$. In the standard coalescent or diffusion limit, let $\tA$
and $\tB$ denote the population-scaled mutation rates at loci $A$ and
$B$, respectively, and let $\rho$ denote the population-scaled
recombination rate between the two loci. Given a sample configuration
$\bfn$ (defined later in the text), consider the following asymptotic
expansion of the sampling probability $q(\bfn\mid\tA,\tB,\gr)$ for
large $\rho$:
%
%
\begin{eqnarray}\label{eqintro,expansion}
q(\bfn\mid\tA,\tB,\gr) &=& q_0(\bfn\mid\tA,\tB) + \frac
{q_1(\bfn\mid\tA,\tB)}{\gr} \nonumber\\[-8pt]\\[-8pt]
&&{} +
\frac{q_2(\bfn\mid\tA,\tB)}{\gr^2} + \cdots,\nonumber
\end{eqnarray}
%
where the coefficients $q_0(\bfn\mid\tA,\tB), q_1(\bfn\mid\tA
,\tB),
q_2(\bfn\mid\tA,\tB),\ldots,$ are independent of $\rho$. The zeroth-order term $q_0$ corresponds to the contribution from the completely
unlinked (i.e., $\rho=\infty$) case, given simply by a product of
marginal one-locus sampling distributions. Until recently, higher-order
terms $q_M(\bfn\mid\tA,\tB)$, for $M\geq1$, were not known. In
\citet{jenson2010}, assuming the infinite-alleles model of
mutation at
each locus, we used probabilistic and combinatorial techniques to
derive a closed-form formula for the first-order term $q_1$, and
showed that the second-order term $q_2$ can be decomposed into two
parts, one for which we obtained a closed-form formula and the other
which satisfies a simple recursion that can be easily evaluated using
dynamic programming. We later extended these results to an arbitrary
finite-alleles model and showed that the same functional form of $q_1$
is shared by all mutation models, a property which we referred to as
\textit{universality} [see \citet{jenson2009G} for details].
Importantly, we also performed an extensive study of the accuracy of
our results and showed that they may be accurate even for moderate
values of $\rho$, including a range that is of biological interest.

Given the above findings, one is naturally led to ask several important
follow-up questions.
In particular, the following questions are of both theoretical and
practical interest.
\begin{enumerate}
\item Is it possible to compute the higher-order coefficients
$q_M(\bfn\mid\tA,\tB)$ for \mbox{$M > 2$}?
\item For a given finite $\rho> 0$, does the series in (\ref
{eqintro,expansion}) converge as more terms are added?
\item If not, how should one make use of the asymptotic expansion in practice?
\item Is it possible to incorporate into the asymptotic expansion
framework other important mechanisms of evolution such as natural
selection?\vadjust{\goodbreak}
\end{enumerate}
In this paper, we develop a new computational technique to answer the
above questions. Our previous method requires rewriting complex
recursion relations into more structured forms, followed by laborious
computation of the expectation of increasingly complicated functions of
multivariate hypergeometric random variables. Generalizing that method
to go beyond the second order (i.e., $M > 2$) seems unwieldy. In
contrast, our new method is based on the diffusion process and it
utilizes the diffusion generator to organize computation in a simple,
transparent fashion. Moreover, the same basic procedure, which is
purely algebraic, applies to all orders and the computation can be
completely automated; we have, in fact, made such an implementation.

To recapitulate, we propose here a method of computing the asymptotic
expansion (\ref{eqintro,expansion}) to an arbitrary order. That is,
for any given positive integer~$M$, our method can be used to compute
the coefficients $q_k(\bfn\mid\tA,\tB)$ for all \mbox{$k \leq M$}; Theorem
\ref{thmg} summarizes this result.
As discussed in Section~\ref{secexample}, however, one can find
examples for which
the series (\ref{eqintro,expansion}) diverges for finite, nonzero
$\rho$. To get around this problem, we employ the method of Pad\'e
approximants. The key idea behind Pad\'e approximants is to approximate
the function of interest by a rational function. Although (\ref
{eqintro,expansion}) may diverge, we show that the sequence of Pad\'e
approximants converges for all values of \mbox{$\rho> 0$}. In fact, for every
sample configuration $\bfn$, we show that there exists a~finite
positive integer $C(\bfn)$, such that the Pad\'e approximant of the
asymptotic expansion up to order $\geq C(\bfn)$ is equal to the \textit
{exact} two-locus sampling distribution. Hence, our result implies that
only a finite number of terms in the asymptotic expansion need to be
computed to recover (via the Pad\'e approximant) the exact sampling
distribution as an analytic function of $\rho$; this function is exact
for all values of $\rho$, including $0$. Theorem \ref{thmconvergence} and
the surrounding discussion lay out the details.
Last, we also show in this paper that our new framework is flexible
enough to incorporate a general model of diploid selection. This
extension is detailed in Section \ref{secselection}.

The above-mentioned convergence result is theoretically appealing. For
practical applications, however, one needs to bear in mind that the
value of $C(\bfn)$ generally grows with sample size, thus implying
that obtaining an exact, analytic sampling distribution may be
impracticable for large samples.
A possible remedy, which works well in practice, is to compute the
asymptotic expansion only up to some reasonable order $M < C(\bfn)$,
and use the corresponding Pad\'e approximant as an approximate sampling
distribution.
We show in Section \ref{secaccuracy} that using $M=10$ or so produces quite
accurate results.

An important advantage of our method over Monte Carlo-based methods is
that, for a given mutation model, the bulk of the computation in our
approach needs to be carried out only once. Specifically, the
coefficients $q_k(\bfn\mid\tA,\tB)$ need to be computed only once,
and the same coefficients can\vadjust{\goodbreak} be used to evaluate the sampling
distribution at different values of the recombination rate $\rho$. We
expect this aspect of our work to have important practical
implications. For example, in the composite likelihood method for
estimating fine-scale recombination rates [\citet{hud2001G},
\citet{mcvetal2004}], one needs to generate exhaustive lookup
tables containing two-locus sampling probabilities for a wide range of
discrete $\rho$ values. An alternative approach would be to store the
coefficients $q_k(\bfn\mid\tA,\tB)$ instead of generating an
exhaustive lookup table using importance sampling, which is
computationally expensive.

The rest of this paper is organized as follows. In Section \ref{secreview},
we lay out the notational convention adopted throughout this paper and
review our previous work on asymptotic expansion of the two-locus
sampling distribution up to second order. Our new technique for
obtaining an arbitrary-order asymptotic expansion is described in
Section \ref{secdiffusion}, where we focus on the selectively neutral
case. In Section~\ref{secOTR,Pade}, we present the method of Pad\'e
approximants and
describe the aforementioned result on convergence to the exact sampling
distribution. In Section \ref{secselection}, we describe how natural
selection can be incorporated into our new framework. Finally, we
summarize in Section \ref{secaccuracy} our empirical study of the accuracy
of various approximate sampling distributions and provide in Section
\ref{secproofs} proofs of the main theoretical results presented in
this paper.

\section{Notation and review of previous work}
\label{secreview}
In this section, we introduce some notation and briefly review previous
results on asymptotic sampling distributions. Initial results were
obtained for the infinite-alleles model of mutation [\citet{jenson2010}]
and later generalized to an arbitrary finite-alleles model
[\citet{jenson2009G}]. In this paper we focus on the latter case.

The set of nonnegative integers is denoted by $\bbN$.
Given a positive integer~$k$, $[k]$~denotes the set $\{1,\ldots,k\}$.
For a~nonnegative real number~$z$ and a~positive integer $n$, $(z)_n
:= z(z+1)\cdots(z+n-1)$ denotes the $n$th ascending factorial of~$z$.
We use $\bfzero$ to denote either a vector or a matrix of all zeroes;
it will be clear from context which is intended. Throughout, we
consider the diffusion limit of a haploid exchangeable model of random
mating with constant population size $2N$. We refer to the haploid
individuals in the population as gametes. Initially we shall assume
that the population is selectively neutral; we extend to the
nonneutral case in Section \ref{secselection}.

\subsection{One-locus sampling distribution}
The sample configuration at a locus is denoted by a vector $\bfn=
(n_1,\ldots,n_K)$, where $n_i$ denotes the number of gametes with
allele $i$ at the locus, and $K$ denotes the number of distinct
possible alleles. We use $n = \sum_{i=1}^K n_i$ to denote the total
sample size. Let $u$ denote the probability of mutation at the locus
per gamete per generation. Then, in the diffusion limit, $N \to\infty
$ and $u \to0$ with the population-scaled mutation rate\vadjust{\goodbreak} $\theta= 4Nu$
held fixed. Mutation events occur according to a~Poisson process with
rate $\theta/2$, and allelic changes are described by a~Markov chain
with transition matrix $\bfP= (P_{ij})$; that is, when a mutation
occurs to an allele $i$, it mutates to allele $j$ with probability $P_{ij}$.

We denote by $p(\bfn)$ the probability of obtaining the \textit
{unordered} sample configuration $\bfn$.
When writing sampling probabilities, we suppress the dependence on
parameters for ease of notation. By exchangeability, the probability of
any \textit{ordered} configuration corresponding to $\bfn$ is invariant
under all permutations of the sampling order. We may, therefore, use
$q(\bfn)$ without ambiguity to denote the probability of any \textit
{particular} ordered configuration consistent with $\bfn$.
The two probabilities are related by
\[
p(\bfn) = \frac{n!}{n_1!\cdots n_K!} q(\bfn).
\]
Throughout this paper, we express our results in terms of ordered
samples for convenience.

Consider an infinite population specified by the population-wide allele
frequencies ${\bfmath{x}}= (x_i,\ldots,x_K)$, evolving according to a
Wright--Fisher diffusion on the simplex
%
%
\begin{equation}\label{eqDeltaK}
\Delta_K = \Biggl\{{\bfmath{x}}=(x_i)\in[0,1]^K \dvtx
\sum_{i=1}^K x_i = 1\Biggr\}.
\end{equation}
We assume that a sample is drawn from the population at stationarity.
No closed-form expression for $q(\bfn)$ is known except in the special
case of parent-independent mutation (PIM), in which $P_{ij}=P_j$ for
all $i$. In the PIM model, the stationary distribution of ${\bfmath
{x}}$ is Dirichlet with parameters $(\theta P_1,\ldots,\theta P_K)$
[\citet{wri1949}], and so $q(\bfn)$ is obtained by drawing an
ordered sample from this population:
%
%
\begin{equation}\label{eqPIM1}
q(\bfn) = \bbE\Biggl[\prod_{i=1}^K X_i^{n_i}\Biggr] = \Gamma
(\theta)\int_{\Delta_K}\prod_{i=1}^K \frac{x_i^{n_i + \theta P_i -
1}}{\Gamma(\theta P_i)} \ud{\bfmath{x}}= \frac{1}{\phammer{\theta
}{n}}\prod_{i=1}^K \phammer{\theta P_i}{n_i}.\hspace*{-30pt}
\end{equation}
This sampling distribution can also be obtained by coalescent arguments
[\citet{gritav1994TPB}].

\subsection{Two loci}
We now extend the above notation to two loci, which we refer to as $A$
and $B$. Denote the probability of a recombination event between the
two loci per gamete per generation by $r$. In the diffusion limit, as
\mbox{$N\to\infty$} we let $r \to0$ such that the population-scaled
recombination parameter $\rho= 4Nr$ is held fixed. Suppose there are
$K$ possible alleles at locus~$A$ and $L$ possible alleles at locus
$B$, with respective population-scaled mutation parameters $\theta_\A
$ and $\theta_\B$, and respective mutation transition matrices $\bfP
^\A$ and~$\bfP^\B$. The two-locus sample configuration is denoted by
$\bfn= (\bfa,\bfb,\bfc)$, where $\bfa= (a_1,\ldots,a_K)$ with\vadjust{\goodbreak}
$a_i$ being the number of gametes with allele $i$ at locus~$A$ and
unspecified alleles at locus $B$; $\bfb= (b_1,\ldots,b_L)$ with $b_j$
being the number of gametes with unspecified alleles at locus $A$ and
allele $j$ at locus~$B$; $\bfc= (c_{ij})$ is a $K \times L$ matrix
with $c_{ij}$ being the multiplicity of gametes with allele $i$ at
locus $A$ and allele $j$ at locus $B$. We also define
\begin{eqnarray*}
a &=& \sum_{i=1}^K a_i,\qquad c_{i\cdot} = \sum_{j=1}^L c_{ij},\qquad
c = \sum_{i=1}^K \sum_{j=1}^L c_{ij}, \\
b &=& \sum_{j=1}^L b_j,\qquad  c_{\cdot j} = \sum_{i=1}^K c_{ij},\qquad
n = a + b + c,
\end{eqnarray*}
and use $\bfc_\A= (c_{i\cdot})_{i\in[K]}$ and $\bfc_\B= (c_{\cdot
j})_{j \in[L]}$ to denote the marginal sample configurations of $\bfc
$ restricted to locus $A$ and locus $B$, respectively. Notice the
distinction between the vectors $\bfa$ and $\bfb$, which represent
gametes with alleles specified at only one of the two loci, and the
vectors $\bfc_\A$ and $\bfc_\B$, which represent the one-locus
marginal configurations of gametes with both alleles observed.

When we consider the ancestry of a sample backward in time, a gamete
may undergo recombination between the two loci, with each of its two
parents transmitting genetic material at only one of the two loci. We
allow the nontransmitting locus to remain unspecified as we trace the
ancestry further back in time.

Denote by $q(\bfa,\bfb,\bfc)$ the sampling probability of an ordered
sample with configuration $(\bfa,\bfb,\bfc)$, again suppressing the
dependence on parameters for ease of notation. Sampling is now from a
two-dimensional Wright--Fisher diffusion with population allele
frequencies ${\bfmath{x}}= (x_{ij})_{i\in[K],j\in[L]}$, evolving on
the state space
%
%
\begin{equation}\label{eqDeltaKL}
\Delta_{K\times L} = \Biggl\{{\bfmath{x}}=(x_{ij}) \in[0,1]^{K\times
L}\dvtx
\sum_{i=1}^K \sum_{j=1}^L x_{ij} = 1\Biggr\}.
\end{equation}
As before, $q(\bfa,\bfb,\bfc)$ is specified by drawing an ordered
sample from the population at stationarity: $q(\bfa,\bfb,\bfc) =
\bbE[F(\bfX;\bfn)]$, where
%
%
\begin{equation}
\label{eqf}
F({\bfmath{x}};\bfn) = \Biggl(\prod_{i=1}^K x_{i\cdot}^{a_i}
\Biggr)\Biggl(\prod_{j=1}^L x_{\cdot j}^{b_j}\Biggr)\Biggl(\prod_{i=1}^K
\prod_{j=1}^L x_{ij}^{c_{ij}}\Biggr)
\end{equation}
with $x_{i\cdot}\!=\!\sum_{j=1}^L x_{ij}$ and $x_{\cdot j}\!=\!\sum
_{i=1}^K x_{ij}$. In the two-locus model with \mbox{$0\!\leq\!\rho\!<\!\infty$},
the stationary distribution, and hence, the sampling distribution, is
not known in closed-form even when the mutation process is
parent-independent. However, when $\rho= \infty$, the two loci become
independent and $q(\bfa,\bfb,\bfc)$ is simply the product of the two
marginal one-locus sampling distributions. More precisely,\vadjust{\goodbreak} denoting the
one-locus sampling distributions at $A$ and $B$ by $q^\A$ and $q^\B$,
respectively, we have
\[
\lim_{\rho\to\infty} q(\bfa,\bfb,\bfc) = q^\A(\bfa+ \bfc_\A
)q^\B(\bfb+ \bfc_\B)
\]
for all mutation models [\citet{eth1979JAP}].
In particular, if mutation is parent-independent, then we do have a
closed-form formula for $q(\bfa,\bfb,\bfc)$ when \mbox{$\rho=\infty$},
since from (\ref{eqPIM1}) we know that
%
%
\begin{equation}
\label{eqmarginalqs}\qquad
q^\A(\bfa) = \frac{1}{\phammer{\theta_\A}{a}}\prod_{i=1}^K
\phammer{\theta_\A P^\A_i}{a_i}  \quad\mbox{and}\quad
q^\B(\bfb) = \frac{1}{\phammer{\theta_\B}{b}}\prod_{j=1}^L
\phammer{\theta_\B P^\B_j}{b_j}.
\end{equation}

\subsection{Asymptotic sampling formula}
\label{secasf}
As mentioned in the \hyperref[int]{Introduction}, although a closed-form formula for
$q(\bfa,\bfb,\bfc)$ is not known for an arbitrary~$\rho$,
previously we [Jenkins and Song (\citeyear{jenson2009G,jenson2010})]
were able to make progress by posing, for large $\rho$, an asymptotic
expansion of the form
%
%
\begin{equation}
\label{eqmainexpansion}
q(\bfa,\bfb,\bfc) = q_0(\bfa,\bfb,\bfc) + \frac{q_1(\bfa,\bfb
,\bfc)}{\rho} + \frac{q_2(\bfa,\bfb,\bfc)}{\rho^2} + \cdots,
\end{equation}
where the coefficients $q_k(\bfa,\bfb,\bfc)$, for all $k \geq0$,
are independent of $\rho$.
We summarize our previous results in the following theorem, specialized
to the case of finite-alleles mutation, which is our interest here.
%
\begin{theorem}[{[\citet{jenson2009G}]}]
\label{thmmain}
In the asymptotic~expan-\break sion~(\ref{eqmainexpansion})~of~the neutral
two-locus sampling formula, the \textup{zero}th-order term is given by
%
%
\begin{equation}\label{eqzerothOrder}
q_0(\bfa,\bfb,\bfc) = q^\A(\bfa+ \bfc_\A)q^\B(\bfb+ \bfc_\B),
\end{equation}
and the first-order term is given by
%
%
\begin{eqnarray}
\label{eqfirstorder}
q_1(\bfa,\bfb,\bfc) & = & \pmatrix{c\cr2}q^\A(\bfa+\bfc_\A)q^\B
(\bfb+\bfc_\B) \nonumber\\
&&{} -q^\B(\bfb+\bfc_\B)\sum_{i=1}^K \pmatrix{c_{i\cdot}\cr2}q^\A
(\bfa+\bfc_\A-\bfe_i) \nonumber\\[-8pt]\\[-8pt]
&&{} -q^\A(\bfa+\bfc_\A)\sum_{j=1}^L \pmatrix{c_{\cdot j}\cr2}q^\B
(\bfb+\bfc_\B-\bfe_j) \nonumber\\
&&{} +\sum_{i=1}^K \sum_{j=1}^L \pmatrix{c_{ij}\cr2}q^\A(\bfa+\bfc_\A
-\bfe_i)q^\B(\bfb+\bfc_\B-\bfe_j),\nonumber
\end{eqnarray}
where $\bfe_i$ is a unit vector with a $1$ at the $i$th entry and
$0$'s elsewhere. Furthermore, the second-order term can be decomposed as
%
%
\begin{equation}
\label{eqsecondorder}
q_2(\bfa,\bfb,\bfc) = \sigma(\bfa,\bfb,\bfc) + q_2(\bfa+\bfc
_\A,\bfb+\bfc_\B,\bfzero),
\end{equation}
where $\sigma(\bfa,\bfb,\bfc)$ is known analytically and $q_2(\bfa
,\bfb,\bfzero)$ satisfies the recursion relation
%
%
\begin{eqnarray} \label{eqq2ab0}
&&
[a(a+\tA-1) + b(b+\tB-1)]q_2(\bfa,\bfb,\bfzero)
\nonumber\\[-2pt]
&&\qquad=\sum_{i=1}^K a_i(a_i-1)q_2(\bfa-\bfe_i,\bfb,\bfzero) \nonumber\\[-2pt]
&&\qquad\quad{}+ \sum
_{j=1}^L b_j(b_j-1)q_2(\bfa,\bfb-\bfe_j,\bfzero) \nonumber\\[-2pt]
&&\qquad\quad{} + \tA\sum_{i=1}^K a_i \sum_{t=1}^K P_{ti}^\A q_2(\bfa-\bfe
_i+\bfe_t,\bfb,\bfzero) \\[-2pt]
&&\qquad\quad{} + \tB\sum_{j=1}^L b_j \sum_{t=1}^L P_{tj}^\B q_2(\bfa,\bfb
-\bfe_j+\bfe_t,\bfzero) \nonumber\\[-2pt]
&&\qquad\quad{} + 4\sum_{i=1}^K \sum_{j=1}^L a_i b_j [ (a-1)(b-1)q^\A(\bfa)q^\B
(\bfb) \nonumber\\[-2pt]
&&\hphantom{\qquad\quad{} + 4\sum_{i=1}^K \sum_{j=1}^L a_i b_j [}
{} - (b-1)(a_i-1)q^\A(\bfa-\bfe
_i)q^\B(\bfb) \nonumber\\[-2pt]
&&\hphantom{\qquad\quad{} + 4\sum_{i=1}^K \sum_{j=1}^L a_i b_j [}
- (a-1)(b_j-1)q^\A(\bfa)q^\B
(\bfb-\bfe_j) \nonumber\\[-2pt]
&&\hphantom{\qquad\quad{} + 4\sum_{i=1}^K \sum_{j=1}^L a_i b_j [}
+ (a_i-1)(b_j-1)q^\A(\bfa-\bfe
_i)q^\B(\bfb-\bfe_j)]\nonumber
\end{eqnarray}
%
with boundary conditions $q_2(\bfe_i,\bfzero,\bfzero) = 0$,
$q_2(\bfzero,\bfe_j,\bfzero) = 0$ and $q_2(\bfe_i,\bfe_j,\bfzero
) = 0$ for all $i\in[K]$ and $j\in[L]$.
\end{theorem}

The expression for $\sigma(\bfa,\bfb,\bfc)$ can be found in
\citet{jenson2009G} and we do not reproduce it here. Notice that
$q_0(\bfa,\bfb,\bfc)$ and $q_1(\bfa,\bfb,\bfc)$ exhibit
\textit{universality}: their dependence on the model of mutation is
subsumed entirely into the one-locus sampling distributions.

The proof of Theorem \ref{thmmain} used coalescent arguments. By
considering the most recent event back in time in the coalescent
process for the sample, it is possible to write down a recursion
relation for the sampling distribution. In a two-locus, finite-alleles
model, the appropriate recursion is a simple modification of the one
introduced by \citet{gol1984} for the \textit{infinite}-alleles model,
also studied in detail by \citet{ethgri1990}. By
substituting~(\ref{eqmainexpansion}) into this recursion, after some lengthy
algebraic manipulation one can obtain the expressions given in
Theorem~\ref{thmmain} [see Jenkins and Song (\citeyear{jenson2009G,jenson2010})
for details].

\section{Arbitrary-order asymptotic expansion}
\label{secdiffusion}
The approach described in Section \ref{secasf} does not generalize
easily.
In what follows, we introduce a new approach\vadjust{\goodbreak} based on the diffusion
approximation. This new method is more transparent and more easily
generalizable than the one used previously, and we illustrate this
point by developing a method for computing the asymptotic expansion
(\ref{eqmainexpansion}) to an arbitrary order.

\subsection{Diffusion approximation of the two-locus model}
Our approach is based on the diffusion process that is dual to the
coalescent process. The generator for the two-locus finite-alleles
diffusion process is
%
%
\begin{eqnarray}\label{eqgenerator1}
\sL &=& \frac{1}{2}\sum_{i=1}^K\sum_{j=1}^L \Biggl[\sum_{k=1}^K\sum
_{l=1}^L x_{ij}(\delta_{ik}\delta_{jl} - x_{kl})\,\pd{x_{kl}} + \tA
\sum_{k=1}^K x_{kj}(P_{ki}^\A- \delta_{ik}) \nonumber\\[-8pt]\\[-8pt]
&&\hspace*{51.7pt}\hphantom{\frac{1}{2}\sum_{i=1}^K\sum_{j=1}^L \Biggl[}{}
+ \tB\sum_{l=1}^L x_{il}(P_{lj}^\B- \delta_{jl}) + \rho
(x_{i\cdot}x_{\cdot j} - x_{ij})\Biggr]\,\pd{x_{ij}},\hspace*{-26pt}
\nonumber
\end{eqnarray}
where $\delta_{ik}$ is the Kronecker delta. For notational
convenience, henceforth, where not specified otherwise, the indices $i$
and $k$ are assumed to take value in $[K]$, while the indices $j$ and
$l$ are assumed to take value in $[L]$.

In what follows, we change to a new set of variables that capture the
decay of dependence between the two loci, an approach originally due to
Ohta and Kimura (\citeyear{ohtkim1969G,ohtkim1969GRC}). Specifically,
the key quantity of
interest is the following definition.
\begin{definition}
The linkage disequilibrium (LD) between allele $i$ at locus $A$ and
allele $j$ at locus $B$ is given by
\[
d_{ij} = x_{ij} - x_{i\cdot} x_{\cdot j} .
\]
\end{definition}

The collection of $(K+1)(L+1) - 1$ new variables is
\[
(x_{1\cdot},\ldots,x_{K\cdot}; x_{\cdot1}, \ldots, x_{\cdot L};
d_{11}, \ldots, d_{\mathit{KL}}).
\]
The diffusion is then constrained to the $(KL-1)$-dimensional simplex
$\Delta_{K \times L}$ by imposing the conditions
%
%
\begin{eqnarray}
\label{eqdagger}
\sum_{i=1}^K x_{i\cdot} &=& 1;\qquad  \sum_{j=1}^L x_{\cdot j} = 1;\qquad
\sum_{i=1}^K d_{ij} = 0\qquad \forall j; \nonumber\\[-8pt]\\[-8pt]
\sum_{j=1}^L d_{ij} &=& 0\qquad \forall i.\nonumber
\end{eqnarray}

\subsection{Rescaling LD}
Since we are interested in the large $\rho$ limit, we expect each
$d_{ij}$ to be small. We introduce the rescaling $\td_{ij} =
\sqrt{\rho} d_{ij}$. The reason for this choice should become clear
from (\ref{eqgenneutral2}) below; in the resulting generator,
there should be a nontrivial interaction between recombination and
genetic drift, that is, they should both act on the fastest timescale.
The leading-order contribution\vadjust{\goodbreak} to the generator, denoted below by
$L^{(2)}$, has contributions from both of these biological processes if
and only if we use this choice for $\td_{ij}$. See \citet{sonson2007}
for another example of this rescaling. By substituting for the new
variables in (\ref{eqgenerator1}) and using (\ref{eqdagger}) for
extensive simplification, the generator can be expressed as
%
%
\begin{equation}
\label{eqgenneutral2}
\widetilde{\sL} = \frac{1}{2}\biggl[\rho L^{(2)} + \sqrt{\rho
}L^{(1)} + L^{(0)} + \frac{1}{\sqrt{\rho}}L^{(-1)}\biggr],
\end{equation}
where the operators in (\ref{eqgenneutral2}) are given by
\begin{eqnarray*}
L^{(2)} & = & \sum_{i,j} \biggl\{ x_{i\cdot}x_{\cdot j}
\biggl[\sum_{k,l} (\gd_{ik}-x_{k\cdot})(\gd_{jl} - x_{\cdot l})\,\pd
{\td_{kl}}\biggr] - \td_{ij}\biggr\}\,\pd{\td_{ij}}, \\
L^{(1)} & = & \sum_{i,j,k,l} [2\td_{ij}(\gd_{ik}-x_{k\cdot
})(\gd_{jl}-x_{\cdot l}) - \gd_{ik}\gd_{jl}\td_{ij} + 2\td
_{il}x_{k\cdot}x_{\cdot j}]\,\pdd{\td_{kl}}{\td_{ij}},
\\
L^{(0)} & = & -\sum_{i,j,k,l} \td_{ij}\td_{kl}\,\pdd{\td_{kl}}{\td
_{ij}} + 2\sum_{i,k,l}[(\gd_{ik}-x_{i\cdot})\td_{kl} - \td
_{il}x_{k\cdot}]\,\pdd{\td_{kl}}{x_{i\cdot}} \\
&&{} +2\sum_{j,k,l}[(\gd_{jl}-x_{\cdot j})\td_{kl} - \td
_{kj}x_{\cdot l}]\,\pdd{\td_{kl}}{x_{\cdot j}} \\
&&{} + \sum_{i,k} x_{i\cdot}(\gd_{ik} - x_{k\cdot})\,\pdd{x_{k\cdot
}}{x_{i\cdot}}+\sum_{j,l} x_{\cdot j}(\gd_{jl} - x_{\cdot l})\,\pdd
{x_{\cdot l}}{x_{\cdot j}} \\
&&{} +\sum_{i,j} \biggl[\theta_\A\sum_k\td_{kj}(P_{ki}^\A- \gd
_{ik}) + \theta_\B\sum_l\td_{il}(P_{lj}^\B- \gd_{jl}) - 2\td
_{ij}\biggr]\,\pd{\td_{ij}}\\
&&{} + \frac{\theta_\A}{2}\sum_{i,k} x_{k\cdot}(P_{ki}^\A- \gd
_{ik})\,\pd{x_{i\cdot}} + \frac{\theta_\B}{2}\sum_{j,l} x_{\cdot
l}(P_{lj}^\B- \gd_{jl})\,\pd{x_{\cdot j}}, \\
L^{(-1)} & = & 2\sum_{i,j}\td_{ij}\,\pdd{x_{i\cdot}}{x_{\cdot j}}.
\end{eqnarray*}
This generator extends that of \citet{ohtkim1969G} from a
$(2\times
2)$- to a $(K\times L)$-allele model. Note that ours differs from that
of \citet{ohtkim1969G} by a factor of two; one unit of time corresponds
to $2N$ (rather than~$N$) generations in our convention.

Recall that our interest is in calculating the expectation at
stationarity of the function $F({\bfmath{x}};\bfn)$ shown in (\ref
{eqf}), which is now viewed as a function of
\[
\widetilde{{\bfmath{x}}}= (x_{1\cdot},\ldots,x_{K\cdot};x_{\cdot
1},\ldots,x_{\cdot L};\td_{11},\ldots,\td_{\mathit{KL}}).
\]
In the same way that the multiplicity matrix $\bfc$ represents
multinomial samples from a population with frequencies $(x_{ij})_{i\in
[K],j\in[L]}$, we introduce an analogous matrix $\bfr= (r_{ij})_{i\in
[K],j\in[L]}$ associated with the variables $(\td_{ij})_{i\in
[K],j\in[L]}$. We further define the marginal vectors $\bfr_\A=
(r_{i\cdot})_{i\in[K]}$ and $\bfr_\B= (r_{\cdot j})_{j\in[L]}$,
where $r_{i\cdot} = \sum_j r_{ij}$ and $r_{\cdot j} = \sum_i
r_{ij}$, analogous to $\bfc_\A$ and $\bfc_\B$. In this notation,
the function $F({\bfmath{x}};\bfn)$ becomes
%
%
\begin{eqnarray}\label{eqfm}
F(\widetilde{{\bfmath{x}}}; \bfn) &=& \Biggl(\prod_{i=1}^K x_{i\cdot
}^{a_i}\Biggr)\Biggl(\prod_{j=1}^L x_{\cdot j}^{b_j}\Biggr)
\Biggl(\prod_{i=1}^K \prod_{j=1}^L \biggl[\frac{\td_{ij}}{\sqrt{\rho}}
+ x_{i\cdot}x_{\cdot j}\biggr]^{c_{ij}}\Biggr) \nonumber\\[-8pt]\\[-8pt]
&=& \sum_{m=0}^c \frac{1}{\rho^{{m/2}}}\!\sum_{\bfr\in
\partition{m}} \biggl[\prod_{i,j}\!\pmatrix{c_{ij}\cr r_{ij}}\biggr]
G^{(m)}(\widetilde{{\bfmath{x}}}; \bfa+\bfc_\A- \bfr_\A,\bfb
+\bfc_\B-\bfr_\B,\bfr),\hspace*{-25pt}\nonumber
\end{eqnarray}
where
%
%
\begin{equation}
\label{eqgm}
G^{(m)}(\widetilde{{\bfmath{x}}}; \bfa,\bfb,\bfr) = \Biggl(\prod
_{i=1}^K x_{i\cdot}^{a_i}\Biggr)\Biggl(\prod_{j=1}^L x_{\cdot
j}^{b_j}\Biggr)\Biggl(\prod_{i=1}^K\prod_{j=1}^L\td_{ij}^{
r_{ij}}\Biggr),
\end{equation}
and the inner summation in (\ref{eqfm}) is over all $K\times L$
matrices $\bfr$ of nonnegative integers whose entries sum to $m$:
\[
\partition{m} = \Biggl\{\bfr\in\bbN^{K\times L}\dvtx\sum
_{i=1}^K\sum_{j=1}^L r_{ij} = m\Biggr\}.
\]
Note that only those matrices which form ``subsamples'' of $\bfc$ have
nonzero coefficient in (\ref{eqfm}); that is, $0 \leq r_{ij} \leq
c_{ij}$ for all $i$ and $j$.

\subsection{The key algorithm}

We now pose an asymptotic\vspace*{1pt} expansion for the expectation $\bbE
[G^{(m)}(\tbfX;\bfa,\bfb,\bfr)]$:
%
%
\begin{eqnarray}
\label{eqexpansion}
\bbE\bigl[G^{(m)}(\tbfX; \bfa,\bfb,\bfr)\bigr] &=& g_0^{(m)}(\bfa,\bfb,\bfr
) + \frac{g_1^{(m)}(\bfa,\bfb,\bfr)}{\sqrt{\rho}}\nonumber\\[-8pt]\\[-8pt]
&&{} + \frac
{g_2^{(m)}(\bfa,\bfb,\bfr)}{\rho} + \cdots,\nonumber
\end{eqnarray}
so that, using (\ref{eqfm}), the quantity of interest is given by
%
%
\begin{eqnarray}\label{eqlevels}
q(\bfa,\bfb,\bfc) &=& \bbE[F(\tbfX; \bfn)] \nonumber\\
&=& \sum_{m=0}^c \sum_{\bfr\in\partition{m}}\Biggl[\prod
_{i=1}^K\prod_{j=1}^L \pmatrix{c_{ij}\cr r_{ij}}\Biggr]\\
&&\hphantom{\sum_{m=0}^c \sum_{\bfr\in\partition{m}}}{}\times\sum_{u
=0}^\infty\frac{g^{(m)}_{u}(\bfa+\bfc_\A-\bfr_\A,\bfb+\bfc
_\B-\bfr_\B,\bfr)}{\rho^{(m+u)/2}}.\nonumber
\end{eqnarray}
We also have the boundary conditions
%
%
\begin{equation}
\label{eqboundaries}
q(\bfe_i,\bfzero,\bfzero) = \pi^\A_i,\qquad  q(\bfzero,\bfe_j,\bfzero
) = \pi^\B_j,\qquad q(\bfe_i,\bfe_j,\bfzero) = \pi^\A_i \pi^\B_j,
\end{equation}
where $\bfmathh{\pi}^\A= (\pi^\A_i)_{i\in[K]}$ and $\bfmathh{\pi
}^\B= (\pi^\B_j)_{j\in[L]}$ are the stationary distributions of
$\bfP^\A$ and $\bfP^\B$, respectively.

Using (\ref{eqlevels}) and (\ref{eqboundaries}), we can
also assign boundary conditions for each $g^{(0)}_u(\bfa,\bfb
,\bfzero)$:
%
%
\begin{eqnarray}
\label{eqnewboundaries}
g^{(0)}_0(\bfe_i,\bfzero,\bfzero) &=& \pi^\A_i,\qquad g^{(0)}_u(\bfe
_i,\bfzero,\bfzero) = 0\qquad \forall u\geq1,\nonumber\\[-2pt]
g^{(0)}_0(\bfzero,\bfe_j,\bfzero) &=& \pi^\B_j,\qquad g^{(0)}_u
(\bfzero,\bfe_j,\bfzero) = 0\qquad \forall u \geq1,\\[-2pt]
g^{(0)}_0(\bfe_i,\bfe_j,\bfzero) &=& \pi^\A_i\pi^\B_j,\qquad
g^{(0)}_u(\bfe_i,\bfe_j,\bfzero) = 0\qquad \forall u\geq
1.\nonumber
\end{eqnarray}

\begin{figure}

\includegraphics{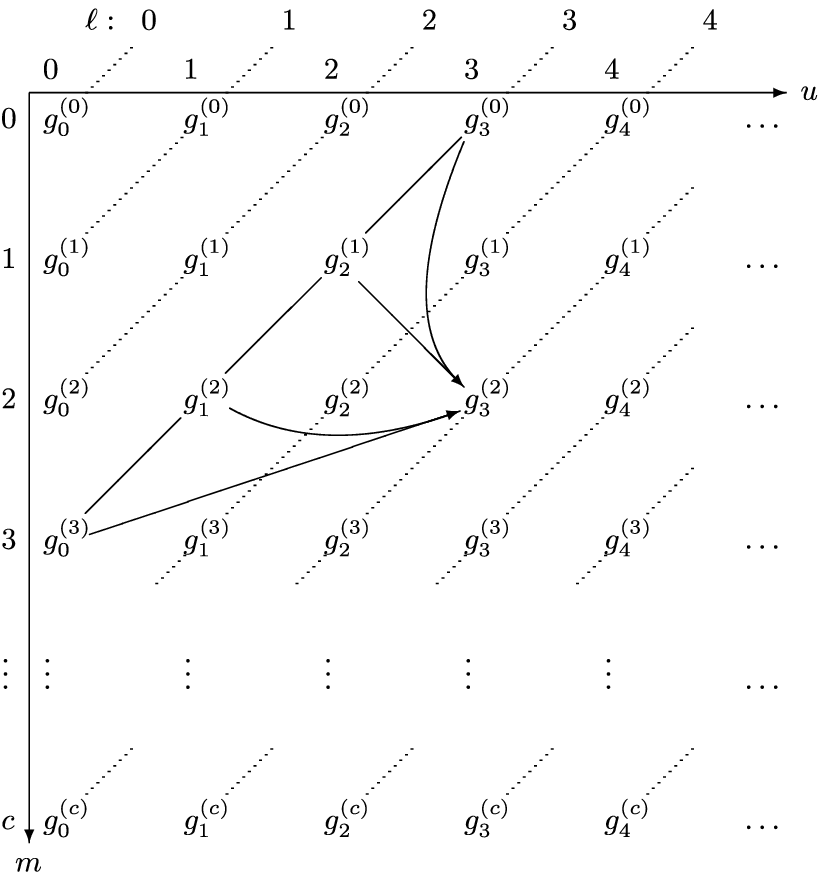}

\caption{Computation of $g_u^{(m)}$ for each
$m=0,\ldots,c$ and $u=0,1,\ldots.$ When $m>0$, the term $g_u ^{(m)}$
residing on level $\ell=m+u$ is determined by up to four entries on
level $\ell-2$. This is illustrated for the example
$g_3^{(2)}$.}\label{figarray}
\vspace*{-3pt}
\end{figure}

We have reduced the problem of computing an asymptotic expansion for
$q(\bfa,\bfb,\bfc)$ to one of computing $g^{(m)}_u(\bfa,\bfb
,\bfr)$ for each $m$ and $u$. Consider arranging these quantities in
a $c \times\bbN$ array, as illustrated in Figure \ref{figarray}.
%
%
Refer to
entries on the $\ell$th anti-diagonal, such that $m+u= \ell$, as
residing on the $\ell$th \textit{level}.
As is clear from (\ref{eqlevels}), the contribution in the expansion
for $q(\bfa,\bfb,\bfc)$ of order $\rho^{-\ell/2}$ is comprised of
entries on level $\ell$. For convenience, we define $g_u^{(m)}(\bfa
,\bfb,\bfr) = 0$ if $ u<0$, or $m<0$, or if any entry
$a_i,b_j,r_{ij}<0$. Then the following theorem, proved
in Section \ref{secproofg}, enables the level-wise computation of each
$g_ u^{(m)}$.\vadjust{\goodbreak}
\begin{theorem}
\label{thmg}
The term $g_ u^{(m)}(\bfa,\bfb,\bfr)$ in the right-hand side of
(\ref{eqexpansion}) is determined as follows.
\begin{longlist}
\item For $m=0$ and $ u=0$, $g_0^{(0)}(\bfa,\bfb,\bfzero) = q^\A
(\bfa)q^\B(\bfb)$. [Recall that $q^\A(\bfa)$ and $q^\B(\bfb)$
are the respective one-locus sampling distributions at locus $A$ and
locus $B$.]
\item For $m > 0$ and $ u\geq0$, $g_ u^{(m)}(\bfa,\bfb,\bfr)$ on
level $\ell$ is determined by at most four entries on level $\ell-2$;
these are $g_ u^{(m-2)}$, $g_{ u-1}^{(m-1)}$, $g_{ u-2}^{(m)}$ and
$g_{ u-3}^{(m+1)}$. This relationship is given explicitly [equation
(\ref{eqcomputeg})] in Section \ref{secproofs}.
\item
\hypertarget{itemcomplication}
For $m=0$ and $ u\geq1$, $g_ u
^{(0)}(\bfa,\bfb,\bfzero)$ on level $\ell$ is determined by $g_{ u
-1}^{(1)}$, also on level $\ell$, and other similar terms $g_ u
^{(0)}(\bfa',\bfb',\bfzero)$, where $a' \leq a$, $b' \leq b$. This
relationship is given explicitly [equation (\ref{eqcomputeg2})] in
Section \ref{secproofs}.
\end{longlist}
\end{theorem}

For odd levels (i.e., $\ell=1,3,5,\ldots$),
the above theorem implies the following vanishing result, which we
prove in Section \ref{secproofodd}.
\begin{corollary}
\label{corodd}
If $m+ u=\ell$ is odd, then $g_ u^{(m)}(\bfa,\bfb,\bfr) = 0$, for
all configurations $(\bfa,\bfb,\bfr)$.
\end{corollary}

Incidentally, Corollary \ref{corodd} implies that only integral powers
of $\rho$ have nonzero coefficients in (\ref{eqlevels}).

Given the entries on level $\ell-2$, Theorem \ref{thmg} provides a method
for computing each of the entries on level $\ell$. They can be
computed in any order, apart from the slight complication of
\hyperlink{itemcomplication}{(iii)}, which requires knowledge of $g_{ u-1}^{(1)}$ as a
prerequisite to $g_ u^{(0)}$. The expression for $g_ u^{(0)}(\bfa
,\bfb,\bfzero)$ is only known recursively in $\bfa$ and~$\bfb$, and
we do not have a closed-form solution for this recursion for $ u\geq
4$ and even. Equation (\ref{eqq2ab0}) provides one example, which
turns out to be the recursion for $g_4^{(0)}(\bfa,\bfb,\bfzero)$.
If the marginal one-locus sampling distributions are known, the
complexity of computing ${g}^{(m)}_{u}(\bfa,\bfb,\bfr)$ for $m > 0$
does not depend on the sample configuration $(\bfa,\bfb,\bfr)$. In
contrast, the complexity of computing $g^{(0)}_{u}(\bfa,\bfb,\bfzero
)$ depends on $\bfa$ and~$\bfb$, and the running time generally grows
with sample size.
However, in Section \ref{secaccuracy} we show that ignoring
$g^{(0)}_{u}(\bfa,\bfb,\bfzero)$ generally leads to little loss of
accuracy in practice.

%
%
\begin{table}
\caption{Entries on levels $\ell=0,2$ of the array $(g_ u^{(m)})$}
\label{tabrecapitulate}
\begin{tabular*}{\tablewidth}{@{\extracolsep{\fill}}lcccc@{}}
\hline
$\bolds\ell$ & $\bolds{u}$ & $\bolds{m}$ & $\bolds{g_ u^{(m)}(\bfa,\bfb,\bfr)}$
& \textbf{Contribution to (\ref{eqmainexpansion})}\\
\hline
$0$ & $0$ & $0$ & $\bbE[\prod_i X_{i\cdot}^{a_i}\prod_j
X_{\cdot j}^{b_j}]$ & $q_0(\bfa,\bfb,\bfc)$ \\
[6pt]
2 & $0$ & $2$ & $\bbE[X_{i\cdot}X_{\cdot
j}(\gd_{ik} - X_{k\cdot})(\gd_{jl} - X_{\cdot l})\prod_u
X_{u\cdot}^{a_u}\prod_v X_{\cdot v}^{b_v}]$, & $q_1(\bfa,\bfb
,\bfc)$\\
& & & \multicolumn{1}{l}{\qquad where $\bfr= \bfe_{ij} + \bfe_{kl}$} & \\
2 & $1$ & $1$ & $0$ & $0$ \\
2 & $2$ & $0$ & $0$ & $0$ \\
\hline
\end{tabular*}
\end{table}
%

To illustrate the method, entries on the first two even levels are
summarized in Table \ref{tabrecapitulate}. These recapitulate part of the
results given in Theorem \ref{thmmain}.
The last column of Table \ref{tabrecapitulate} gives the ``contribution''
to (\ref{eqmainexpansion}) from $g_ u^{(m)}(\bfa,\bfb,\bfr)$ (for
fixed $ u$ and $m$ and summing over the relevant $\bfr$). According
to (\ref{eqlevels}), this quantity is
%
%
\begin{equation}
\label{eqmsum}
\sum_{\bfr\in\partition{m}}\Biggl[\prod_{i=1}^K\prod_{j=1}^L
\pmatrix{c_{ij}\cr r_{ij}}\Biggr] g_ u^{(m)}(\bfa+\bfc_\A-\bfr_\A
,\bfb+\bfc_\B-\bfr_\B,\bfr).
\end{equation}
We have also checked that the total contribution from entries on level
$\ell= 4$ is equal to $q_2(\bfa,\bfb,\bfc)$, as given in Theorem
\ref{thmmain}. We note in passing that Theorem \ref{thmg} makes it
transparent why $q_2(\bfa,\bfb,\bfc)$ is \textit{not} universal in
the sense described in Section \ref{secasf}: expressions on level
$\ell= 4$
depend directly on $L^{(0)}$, which in turn depends upon the model of
mutation. By contrast, the nonzero contribution to $q_1(\bfa,\bfb
,\bfc)$, for example, is determined by $L^{(2)}$, which does not
depend on the model of mutation.

It is important to emphasize that $g_ u^{(m)}(\bfa,\bfb,\bfr)$ is a
function of the vectors $\bfa$, $\bfb$, the matrix $\bfr$ and
(implicitly) the parameters $\theta_\A$ and $\theta_\B$. The
relationships given in Theorem \ref{thmg} are thus \textit
{functional}, and
only need to be computed once. In other words, all of these arguments
of $g_ u^{(m)}$ can remain arbitrary.
It is not necessary to redo any of the algebraic computations for each
particular choice of sample configuration, for example. Moreover, the
solutions to each $g_ u^{(m)}(\bfa,\bfb,\bfr)$ are expressed
concisely in terms of the marginal one-locus sampling distributions
$q^\A$ and $q^\B$; this fact follows inductively from the solution
for $g_0^{(0)}(\bfa,\bfb,\bfzero)$.
Unlike the method of Jenkins and Song (\citeyear{jenson2009G,jenson2010}),
the iterative procedure here is essentially the same at every step.

\section{\texorpdfstring{Partial sums, optimal truncation and Pad\'e
approximants}{Partial sums, optimal truncation and Pade approximants}}
\label{secOTR,Pade}

In principle, the procedure described in the previous section provides
a method of computing an arbitrary number of terms in the asymptotic
expansion (\ref{eqmainexpansion}), for any sample configuration.
Suppose the computation has been carried out up to level $\ell=2M$ and consider
the partial sum
%
%
\begin{equation}
\label{eqpartialsum}
\qPS{M}(\bfa,\bfb,\bfc) = q_0(\bfa,\bfb,\bfc) + \frac{q_1(\bfa
,\bfb,\bfc)}{\rho} + \cdots+ \frac{q_{M}(\bfa,\bfb,\bfc)}{\rho^{M}}.
\end{equation}
%
Since we do not know its radius of convergence, we should be prepared
for this sum to diverge eventually as $M$ increases.
An important question that we address in this section is: How many
terms should we use to maximize the accuracy of the approximation?

\subsection{Optimal truncation}
\label{secOTR}
As mentioned above, simply adding more and more terms to the asymptotic
expansion may decrease the accuracy beyond a certain point.
Optimal truncation is a rule of thumb for truncating the partial sum at
a point that is expected to maximize its accuracy.
More precisely, it is defined as follows.
\begin{definition}[(Optimal truncation rule)] Given the first $M+1$ terms
$q_0(\bfa,\bfb,\bfc)$, $q_1(\bfa,\bfb,\bfc), \ldots,
q_{M}(\bfa,\bfb,\bfc)$, in the asymptotic expansion,\vspace*{1pt} let $M'$ be
the index such that $| {q_{M'}(\bfa,\bfb,\bfc)}/{\rho
^{M'}}| < | {q_{M''}(\bfa,\bfb,\bfc)}/{\rho
^{M''}}|$, for all \mbox{$M'' \neq M'$}, where $M', M'' \leq M$. Then,
the optimal truncation rule (OTR) suggests truncating the sum at
order $M'$:
\[
\qOTR{M}(\bfa,\bfb,\bfc) = q_0(\bfa,\bfb,\bfc) + \frac{q_1(\bfa
,\bfb,\bfc)}{\rho} + \cdots+ \frac{q_{M'}(\bfa,\bfb,\bfc)}{\rho^{M'}}.
\]
\end{definition}

The motivation for this rule is that the magnitude of the
$i$th term in the expansion is an estimate of the magnitude of the
remainder. More sophisticated versions of the OTR are available
[e.g., \citet{din1973}, Chapter XXI], but for simplicity we focus
on the
definition given above.

There are two issues with the OTR. First, it minimizes the error only
approximately, and so, despite its name, it is not guaranteed to be
optimal. For example, the magnitude of the first few terms in the
series can behave very irregularly before a pattern emerges. Second, it
may use only the first few terms in the expansion and discard the rest.
As we discuss later, for some sample configurations and parameter
values of interest, the OTR might truncate very early, even as early as
$M'=2$. This is unrelated to the first issue, since the series may
indeed begin to diverge very early.
Below, we discuss a better approximation scheme with a provable
convergence property.

\subsection{\texorpdfstring{Pad\'e approximants}{Pade approximants}}
\label{secPade}
The key idea behind Pad\'e approximants is to approximate the function
of interest by a rational function. In contrast to the OTR, Pad\'e
approximants make use of \textit{all} of the computed terms in the
expansion, even when the expansion diverges rapidly. More precisely,
the $[U/V]$ Pad\'e approximant of a function is defined as follows.
\begin{definition}[({$[U/V]$} Pad\'e approximant)]
Given a function $f$ and two nonnegative integers $U$ and $V$, the
$[U/V]$ Pad\'e approximant of $f$ is a~rational function of the form
\[
\pade{U/V}{f}{x} = \frac{A_0 + A_1 x + \cdots+ A_U x^U}{B_0 + B_1 x
+ \cdots+ B_V x^V},
\]
such that $B_0=1$ and
\[
f(x) - \pade{U/V}{f}{x} = O(x^{U+V+1}).
\]
That is, the first $U+V+1$ terms in a Maclaurin series of the Pad\'e
approxi\-mant $\pade{U/V}{f}{x}$ matches the first $U\!+\!V\!+\!1$ terms in a
Maclaurin series of~$f$.\vadjust{\goodbreak}
\end{definition}

Our goal is to approximate the sampling distribution $q(\bfa,\bfb
,\bfc)$ by the Pad\'e approximant
%
%
\begin{equation}\label{eqqPade}
\qPade{U/V}(\bfa,\bfb,\bfc) = [U/V]_{q(\bfa,\bfb,\bfc)}
\biggl(\frac{1}{\rho}\biggr),
\end{equation}
such that the first $U+V+1$ terms in a Maclaurin series of $\pade
{U/V}{q(\bfa,\bfb,\bfc)}{\frac{1}{\rho}}$ agrees with
(\ref{eqpartialsum}), where $M = U+V$.
[In this notation, $\pade{U/V}{q(\bfa,\bfb,\bfc)}{\frac{1}{\rho
}}$ is an implicit function of the mutation parameters.]
As more terms in (\ref{eqpartialsum}) are computed (i.e., as $M$
increases), a sequence of Pad\'e approximants can be constructed. This
sequence often has much better convergence properties than (\ref
{eqpartialsum}) itself [\citet{bakgra1996}].

For a given $M$, there is still some freedom over the choice of $U$ and~$V$.
As $M$ increases, we construct the following ``staircase''
sequence of Pad\'e approximants: $[0/0], [0/1], [1/1], [1/2],
[2/2],\ldots.$ This scheme is motivated by the following lemma, proved
in Section \ref{secproofrational}.
\begin{lemma}
\label{lemrational}
Under a neutral, finite-alleles model, the sampling distribution
$q(\bfa,\bfb,\bfc)$ is a rational function of $1/\rho$, and the
degree of the numerator is equal to the degree of the denominator.
\end{lemma}

This simple yet powerful observation immediately leads to a convergence
result for the Pad\'e approximants in the following manner.
\begin{theorem}
\label{thmconvergence}
Consider a neutral, finite-alleles model. For every given two-locus
sample configuration $(\bfa,\bfb,\bfc)$,
there exists a finite nonnegative integer $U_0$ such that for all $U
\geq U_0$ and $V \geq U_0$, the Pad\'e approximant $\qPade{U/V}(\bfa
,\bfb,\bfc)$
is exactly equal to $q(\bfa,\bfb,\bfc)$ for all $\rho\geq0$.
\end{theorem}

A proof of this theorem is provided in
Section \ref{secproofconvergence}. Note that the staircase sequence
is the
``quickest'' to reach the true sampling distribution
$q(\bfa,\bfb,\bfc)$. Although Theorem \ref{thmconvergence}
provides very
strong information about the convergence of the Pad\'e approximants, in
practice $U$ and $V$ will have to be intractably large for such
convergence to take place. The real value of Pad\'e summation derives
from the empirical observation that the approximants exhibit high
accuracy even \textit{before} they hit the true sampling distribution.
The staircase sequence also has the advantage that it exhibits a
continued fraction representation, which enables their construction to
be made computationally more efficient [\citet{bakgra1996},
Chapter 4].

\section{Incorporating selection}
\label{secselection}
We now incorporate a model of diploid selection into the results of
Section \ref{secdiffusion}. Suppose that a diploid individual is
composed of
two haplotypes $(i,j)$, $(k,l) \in[K]\times[L]$, and that, without
loss of generality,\vadjust{\goodbreak} selective differences exist at locus $A$. We denote
the fitness of this individual by $1+s^\A_{ik}$, and consider the
diffusion limit in which $\sigma^\A_{ik} = 4Ns^\A_{ik}$ is held
fixed for each $i,k\in[K]$ as $N\to\infty$.

In what follows, we use a subscript ``s'' to denote selectively
nonneutral versions of the quantities defined above. Results will be
given in terms of the nonneutral one-locus sampling distribution $q^\A
_\ts$ at locus $A$ and the neutral one-locus sampling distribution
$q^\B$ at locus $B$.

\subsection{One-locus sampling distribution under selection}
For the infinite-alleles model, one-locus sampling distributions under
symmetric selection have been studied by \citet{gro2002app},
\citet{han2005sam} and \citet{hui2007ewe}. In the case of a
parent-independent finite-alleles model, the stationary distribution of
the one-locus selection model is known to be a \textit{weighted}
Dirichlet distribution [\citet{wri1949}],
\[
\pi^\A_\ts({\bfmath{x}}) = D\Biggl(\prod_{i=1}^K x_i^{\theta_\A
P^\A_i - 1} \Biggr)\exp\Biggl(\frac{1}{2}\sum_{i=1}^K \sum
_{k=1}^K \gs^\A_{ik} x_ix_k\Biggr),
\]
where ${\bfmath{x}}\in\Delta_K$ [see (\ref{eqDeltaK})] and $D$ is a
normalizing constant. The one-locus sampling distribution at
stationarity is then obtained by drawing a multinomial sample from this
distribution:
%
%
\begin{equation}
\label{eqq-s}
q^\A_\ts(\bfa) = \bbE_{\ts}\Biggl[\prod_{i=1}^K X_i^{a_i}\Biggr].
\end{equation}
Thus, under a diploid selection model with parent-independent mutation,
we are able to express the one-locus sampling distribution at least in
integral form. There are two integrals that need to be evaluated: one
for the expectation and the other for the normalizing constant. In
practice, these integrals must be evaluated using numerical methods
[\citet{buzetal2009}].

\subsection{Two-locus sampling distribution with one locus under selection}
To incorporate selection into our framework, we first introduce some
further notation.
\begin{definition}
Given two-locus population-wide allele frequencies ${\bfmath{x}}\in
\Delta_{K\times L}$ [see (\ref{eqDeltaKL})],
the mean fitness of the population at locus $A$ is
\[
\bar{\gs}^\A({\bfmath{x}}_\A) = \sum_{i,k}\gs^\A_{ik}
x_{i\cdot} x_{k\cdot}.
\]
\end{definition}

Selection has an additive effect on the generator of the process, which
is now given by
\[
\sL_{\ts} = \sL+\frac{1}{2}\sum_{i,j}x_{ij}\biggl[\sum_k \gs^\A
_{ik}x_{k\cdot} - \bar{\gs}^\A({\bfmath{x}}_\A)\biggr]\,\pd{x_{ij}},
\]
where $\sL$ is the generator (\ref{eqgenerator1}) of the neutral
diffusion process [see, e.g., \citet{ethnag1989}].
Rewriting $\sL_\ts$ in terms of the LD variables and then rescaling
$d_{ij}$ as before, we obtain
\[
\widetilde{\sL}_\ts= \widetilde{\sL} + \frac{1}{2}\biggl[\rho
L_\ts^{(2)} + \sqrt{\rho}L_\ts^{(1)} + L_\ts^{(0)} + \frac
{1}{\sqrt{\rho}}L_\ts^{(-1)}\biggr],
\]
where $\widetilde{\sL}$ is as in (\ref{eqgenneutral2}), and the new
contributions are
\begin{eqnarray*}
L_\ts^{(2)} &=& 0,\\
L_\ts^{(1)} &=& 0,\\
L_\ts^{(0)} &=& \sum_{i,j}\biggl[d_{ij}\biggl(\sum_k \gs_{ik}^\A
x_{k\cdot} - \bar{\gs}^\A({\bfmath{x}}_\A)\biggr) - \sum
_{k,k'}d_{kj}\gs_{kk'}^\A x_{i\cdot}x_{k'\cdot}\biggr]\,\pd{d_{ij}}\\
&&{} + \sum_i x_{i\cdot}\biggl(\sum_k \gs_{ik}^\A x_{k\cdot} - \bar
{\gs}^\A({\bfmath{x}}_\A)\biggr)\,\pd{x_{i\cdot}},\\
L_\ts^{(-1)} &=& \sum_{i,j,k}\td_{ij} \gs_{ik}^\A x_{k\cdot}\,\pd
{x_{\cdot j}}.
\end{eqnarray*}
In addition, we replace the asymptotic expansion (\ref{eqexpansion}) with
\begin{eqnarray*}
\bbE_\ts\bigl[G^{(m)}(\tbfX; \bfa,\bfb,\bfr)\bigr] &=& h_0^{(m)}(\bfa,\bfb
,\bfr) + \frac{h_1^{(m)}(\bfa,\bfb,\bfr)}{\sqrt{\rho}} \\
&&{}+ \frac
{h_2^{(m)}(\bfa,\bfb,\bfr)}{\rho} + \cdots,
\end{eqnarray*}
the corresponding expansion for the expectation of $G^{(m)}(\tbfX;
\bfa,\bfb,\bfr)$ with respect to the stationary distribution
\textit{under selection} at locus $A$. Finally, the boundary conditions
(\ref{eqnewboundaries}) become the following [\citet{fea2003JAP}]:
%
%
\begin{eqnarray}
\label{eqnewboundaries-s}
h^{(0)}_0(\bfe_i,\bfzero,\bfzero) &=& \phi^\A_i,\qquad h^{(0)}_ u(\bfe
_i,\bfzero,\bfzero) = 0\qquad \forall u\geq1,\nonumber\\
h^{(0)}_0(\bfzero,\bfe_j,\bfzero) &=& \pi^\B_j,\qquad h^{(0)}_ u
(\bfzero,\bfe_j,\bfzero) = 0\qquad \forall u\geq1,\\
h^{(0)}_0(\bfe_i,\bfe_j,\bfzero) &=& \phi^\A_i\pi^\B_j,\qquad
h^{(0)}_ u(\bfe_i,\bfe_j,\bfzero) = 0\qquad \forall u\geq
1,\nonumber
\end{eqnarray}
where $\bfmathh{\phi}^\A= (\phi^\A_i)_{i\in[K]}$ is the stationary
distribution for drawing a sample of size one from a single selected locus.

With only minor modifications to the arguments of Section \ref
{secdiffusion}, each term in the array for $h_ u^{(m)}$ can be
computed in a manner similar to Theorem \ref{thmg}.
In particular, entries on odd levels are still zero. Furthermore, as
proved in Section \ref{secproofmain-s}, we can update Theorem \ref
{thmmain} as follows.
\begin{theorem}
\label{thmmain-s}
Suppose locus $A$ is under selection, while locus $B$ is selectively neutral.
Then, in the asymptotic expansion (\ref{eqmainexpansion}) of the
two-locus sampling distribution, the \textup{zero}th- and first-order
terms are given by~(\ref{eqzerothOrder}) and (\ref{eqfirstorder}),
respectively, with $q^\A_\ts(\bfa)$ in place of
$q^\A(\bfa)$. Furthermore, the second-order term (\ref{eqsecondorder})
may now be decomposed into two parts:
%
%
\begin{equation}
\label{eqsecondorder-s}
q_{2,\ts}(\bfa,\bfb,\bfc) =
q_{2,\ts}(\bfa+\bfc_\A,\bfb+\bfc_\B,\bfzero) +
\sigma_\ts(\bfa,\bfb,\bfc),
\end{equation}
where $\sigma_\ts(\bfa,\bfb,\bfc)$ is given by a known analytic
expression and $q_{2,\ts}(\bfa,\bfb,\bfzero)$ satisfies a slightly
modified version of the recursion relation (\ref{eqq2ab0}) for
$q_2(\bfa,\bfb,\bfzero)$. (These expressions are omitted for brevity.)
\end{theorem}

We remark that the above arguments can be modified to allow for
locus~$B$ also to be under selection, provided the selection is
independent,\vspace*{1pt} with no epistatic interactions, and provided one can
substitute $\phi_j^\B$ for $\pi_j^\B$ in (\ref{eqnewboundaries-s}).
Then, one\vspace*{-1pt} could also simply substitute $q^\B_\ts(\bfb)$ for
$q^\B(\bfb)$ in the expressions for $q_0(\bfa,\bfb,\bfc)$ and
$q_1(\bfa,\bfb,\bfc)$. However, for the boundary conditions
(\ref{eqnewboundaries-s}) to be modified in this way we would need to
extend the result of \citet{fea2003JAP} to deal with \textit{two}
nonneutral loci, and we are unaware of such a result in the
literature.

\section{Empirical study of accuracy}
\label{secaccuracy}
In this section, we study empirically the accuracy of the approximate
sampling distributions discussed in Section \ref{secOTR,Pade}.

\subsection{Computational details}
\label{secimplementation}
As discussed earlier, a major advantage of our technique is that, given
the first $M$ terms in the asymptotic expansion~(\ref{eqmainexpansion}),
the $(M+1)$th term can be found and has to be computed only once. There
are two complications to this statement: first, as mentioned in the
discussion following Theorem \ref{thmg}, the $M$th-order term $q_M$ for
$M\geq2$ has a contribution [namely, $g^{(0)}_{2M}(\bfa,\bfb,\bfzero
)$] that is \textit{not} known in closed form, and is only given recursively.
[Recall that the $M=1$ case is an exception, with ${g}^{(0)}_{2}(\bfa
,\bfb,\bfzero)=0$ for all $(\bfa,\bfb,\bfzero)$.]
In \citet{jenson2009G} it was observed that the contribution of
${g}^{(0)}_{4}(\bfa,\bfb,\bfzero)$ to $q_2(\bfa,\bfb,\bfc)$ is
generally very small, but that its burden in computational time
increases with sample size. Extrapolating this observation to higher-order terms,
we consider making the following approximation.
\begin{approximation}
\label{approxrecursionapproximation}
For all $M \geq2$, assume
\[
{g}^{(0)}_{2M}(\bfa,\bfb,\bfzero) \approx0
\]
for all configurations $(\bfa,\bfb,\bfzero)$.
\end{approximation}

As we show presently, adopting this approximation has little effect on
the accuracy of asymptotic sampling distributions. In what follows, we
use the symbol\vadjust{\goodbreak} ``$\mathring{\phantom{q}}$'' to indicate when the above
approximation is employed. For example, the partial sum $\qPS{M}(\bfa
,\bfb,\bfc)$ in (\ref{eqpartialsum}) becomes $\qaPS{M}(\bfa,\bfb
,\bfc)$ under Approximation \ref{approxrecursionapproximation}.

Upon making the above approximation, it is then possible to construct a
closed-form expression for each subsequent term $\qApprox_M(\bfa,\bfb
,\bfc)$. However, there is a second issue: Given the effort required
to reach the complicated expression for $\qApprox_2(\bfa,\bfb,\bfc
)$ [\citet{jenson2009G}], performing the computation by hand for
$M > 2$ does not seem tractable. Symbolic computation using computer
software such as \textit{Mathematica} is a feasible option, but we defer
this for future work. Here, we are interested in comparing the accuracy
of asymptotic sampling distributions with the true likelihood.
Therefore, we have implemented an automated \textit{numerical}
computation of each subsequent term in the asymptotic expansion, for a
given fixed sample configuration and fixed mutation parameters. For the
samples investigated below, this did not impose undue computational
burden, even when repeating this procedure across \textit{all} samples
of a given size.
Exact numerical computation of the true likelihood is possible for only
small sample sizes (say, up to thirty), so we restrict our study to
those cases.

For simplicity, we assume in our empirical study that all alleles are
selectively neutral. Furthermore, we assume a symmetric, PIM model so
that $\bfP^\A= \bfP^\B= \left(
{1/2\atop 1/2} \enskip{1/2 \atop 1/2}
\right)$ and take $\tA= \tB= 0.01$. This is a
reasonable approximation for modeling neutral single nucleotide
polymorphism (SNP) data in humans [e.g., \citet{mcvetal2002}].
For the
PIM model, recall that the marginal one-locus sampling distributions
are available in closed-form, as shown in
(\ref{eqmarginalqs}).

\subsection{Rate of convergence: An example}
\label{secexample}
\begin{figure}

\includegraphics{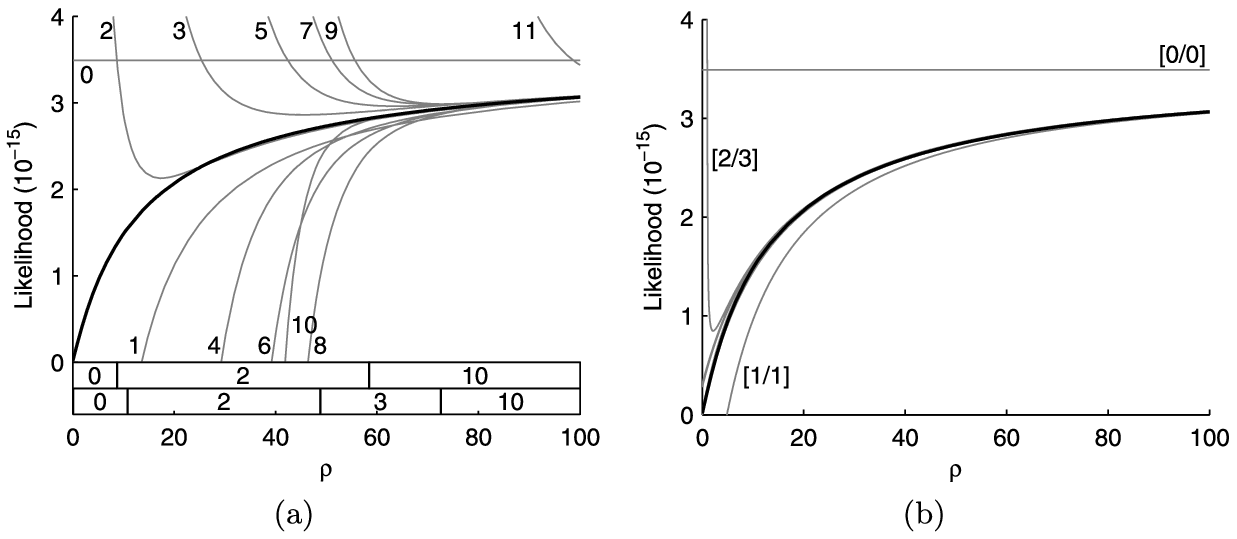}

\caption{Likelihood curves for $\rho$, comparing different levels of
truncation, for the sample $\bfa=\bfb=\bfzero$, $\bfc=
\left({10\atop 2}\enskip{7\atop 1}\right)$.
A symmetric PIM model with $\tA= \tB= 0.01$ is assumed.
The true likelihood is shown as a thick black line. Approximate
likelihood curves for various levels of truncation, $M = 0,1,\ldots
,11$ are shown as thinner gray lines.
\textup{(a)} Partial sums $\qPS{M}$, labeled by $M$. The bottom\vspace*{-1pt} row of
indices records $M'$, the level of truncation in $\qOTR{11}(\abc)$
recommended by the OTR. The row above records the actual level of
truncation that minimizes the unsigned relative error.
\textup{(b)} Staircase Pad\'e approximants $\qPade{U/V}$, some labeled by $[U/V]$.
To the naked eye, $\qPade{U/V}$ for most $[U/V]$ are indistinguishable
from the true likelihood curve.}
\label{figexample}
\end{figure}

To compare the convergence of the sequence of partial sums
(\ref{eqpartialsum}) with that of the sequence of Pad\'e
approximants~(\ref{eqqPade}), we re-examine in detail an example studied previously.
The sample configuration for this example is $\bfa= \bfb= \bfzero$,
$\bfc= \left({10\atop 2}\enskip{7\atop 1}\right)$.
In \citet{jenson2009G}, we were able to
compute the first three terms in the asymptotic expansion, obtaining
the partial sum $\qPS{2}(\abc)$. Applying the new method described in
this paper, we computed $\qPS{M}(\abc)$ for $M \leq11$ [including the
recursive terms $g^{(0)}_{u}(\bfa,\bfb,\bfzero)$ discussed above], and
also the corresponding staircase Pad\'e approximants
$\qPade{U/V}(\abc)$. Results are illustrated in Figure \ref
{figexample}, in
which we compare various approximations of the likelihood curve for
$\rho$ with the true likelihood curve. Here, the likelihood of a
sample is defined simply as its sampling distribution
$q(\bfa,\bfb,\bfc)$ treated as a function of $\rho$, with $\tA$ and
$\tB$ fixed at $0.01$.

%
%
%

Figure \ref{figexample}(a) exhibits a number of features which we also
observed more generally for many other samples. For fixed $\rho$ in
the range illustrated, the sequence $(\qPS{0},\qPS{1},\qPS{2},\ldots
)$ of partial sums diverges eventually, and, for many realistic
choices\vadjust{\goodbreak}
of $\rho$, this divergence can be surprisingly early. Arguably the
best overall fit in the figure is $\qPS{2}$, though for $\rho\geq70$
good accuracy is maintained by $\qPS{M}$ for all $1\leq M \leq10$.
Divergence is extremely rapid if we add any further terms; witness the
curve for $\qPS{11}$. Of course, in real applications the true
likelihood will be unavailable and we might rely on the aforementioned
optimal truncation rule to guide the truncation point. Here, it
performs reasonably well, correctly picking the most accurate index
across most of the range [compare the bottom two rows of indices in
Figure \ref{figexample}(a)].

In contrast, Figure \ref{figexample}(b) shows that there is much less risk
in using ``too many'' terms in constructing the Pad\'e approximants.
They approximate the true likelihood very accurately and very quickly;
to the naked eye, Pad\'e approximants $\qPade{U/V}$ for most $[U/V]$
are indistinguishable from the true likelihood curve. Indeed, this
example suggests that there is very little gain in accuracy beyond
$\qPade{0/1}$, but that there is no significant loss beyond it either.
(It should be pointed out, however, that achievement of such high
accuracy so early in the sequence of approximants does not seem to be a
common occurrence across samples. Usually, several more terms are
required; see next section for further details.)

%
%
\begin{table}
\caption{Nonnegative, real roots in the numerator and denominator of
the $[U/V]$ Pad\'e approximants for the sample $\bfa= \bfb= \bfzero
$, $\bfc= \left({10\atop 2}\enskip{7\atop 1}\right)$}
\label{tabroots}
\begin{tabular*}{\tablewidth}{@{\extracolsep{\fill}}ld{2.3}rrr@{}}
\hline
$\bolds{[U/V]}$ & \multicolumn{2}{c}{\textbf{Roots of numerator}} & \multicolumn
{2}{c@{}}{\textbf{Roots of denominator}}\\
\hline
$[0/0]$ \\
$[0/1]$ & 0.000 \\
$[1/1]$ & 4.871 \\
$[1/2]$ & 0.000 \\
$[2/2]$ \\
$[2/3]$ & 0.000 && 0.912\\
$[3/3]$ \\
$[3/4]$ & 0.000 \\
$[4/4]$ \\
$[4/5]$ & 0.000 \\
$[5/5]$ \\
$[5/6]$ & 0.000 & 8,474,538.140 & 8,474,538.140 \\
$[6/6]$ \\
$[6/7]$ & 0.000 & 306.846 & 306.846\\
$[7/7]$ \\
$[7/8]$ & 0.000 \\
$[8/8]$ & 82.033 && 82.032\\
$[8/9]$ & 0.000 & 77.366 & 0.284 &
77.364\\
$[9/9]$ & 82.121 & 4,412.751 & 82.120 &
4,412.751\\
$[9/10]$ & 0.000 \\
$[10/10]$ & 5.543 && 4.252\\
\hline
\end{tabular*}
\end{table}
%

There is one further important observation to be made regarding Pad\'e
approximants. There is nothing to prevent the polynomial in the
denominator from having positive\vadjust{\goodbreak} real roots, thus creating
singularities. This is indeed observed once in Figure \ref{figexample}(b);
$\qPade{2/3}$ exhibits a singularity at $\rho= 0.9$. To examine this
behavior further, in Table \ref{tabroots} we tabulate the
nonnegative real
roots of the numerator and denominator of each approximant in the
staircase sequence up to $[10/10]$. We see some interesting patterns:
(i)~The total number of nonnegative real roots does not seem to grow
with $M$. (ii)~Roots in the denominator are almost invariably
accompanied by a nearby root in the numerator, and their proximity is
usually extremely close, with agreement to several decimal places.
(iii) Pairs of roots appear transiently; the roots of one approximant
in the sequence provide almost no information regarding the next. Such
pairs of roots are known as \textit{defects}, and are an inevitable
feature of Pad\'e approximants [\citet{bakgra1996}, page 48]. Provided
we anticipate them, defects need not act as a serious obstacle. Because
of the proximity of the two roots in a pair, they almost behave like
removable singularities. In particular, Pad\'e approximants are still
well behaved outside a reasonably narrow interval around the defect,
and can still approximate the likelihood accurately. In light of these
observations, henceforth we use the following simple heuristic for
dealing with defects.\vadjust{\goodbreak}
\begin{definition}[(Defect heuristic)]
\label{defdefect}
Suppose we are interested in approximating the likelihood curve at a
particular value $\rho_0$ and have the resources to compute up to $M$
in the partial sum (\ref{eqpartialsum}). Then proceed as follows.
\begin{longlist}[(3)]
\item[(1)] Initialize with $M'' = M$.
\item[(2)]\hypertarget{step2}
Construct the $[U/V]$ Pad\'e approximant in the
staircase sequence, where $U + V = M''$.
\item[(3)] If it exhibits a root in the interval $(\rho_0 - \varepsilon, \rho
_0 + \varepsilon) \cap[0,\infty)$, either in its numerator or
denominator, then decrement $M''$ by one and go to step~\hyperlink{step2}{(2)},
otherwise use this approximant.
\end{longlist}
\end{definition}

Choice of threshold $\varepsilon$ involves a trade-off between the
disruption caused by the nearby defect and the loss incurred by
reverting to the previous Pad\'e approximant. Throughout the following
section, we employ this heuristic with $\varepsilon= 25$, which seemed to
work well.

\subsection{Rate of convergence: Empirical study}
\label{secsimulation}
To investigate to what extent the observations of Section \ref{secexample}
hold across all samples, we performed the following empirical study.
Following \citet{jenson2009G}, we focused on samples of the form
$(\bfzero,\bfzero,\bfc)$ for which all alleles are observed at both
loci, and measured the accuracy of the partial sum (\ref{eqpartialsum})
by the \textit{unsigned relative error}
\[
\ePS{M}(\bfzero,\bfzero,\bfc) = \biggl| \frac{\qPS{M}(\bfzero
,\bfzero,\bfc) - q(\bfzero,\bfzero,\bfc)}{q(\bfzero,\bfzero,\bfc
)}\biggr| \times100\%.
\]
An analogous definition can be made for $e_{\Pade}^{(M)}$, the
unsigned relative error of the staircase Pad\'e approximants $\qPade
{U/V}$, where $U=\floor{M/2}$ and $V=\ceil{M/2}$. When
Approximation\vspace*{1pt}
\ref{approxrecursionapproximation} is additionally used, the respective
unsigned errors are denoted $\eaPS{M}$ and $\eaPade{M}$. These
quantities are implicit functions of the parameters and of the sample
configuration. In our study, we focused on sufficiently small sample
sizes so that the true sampling distribution $q(\bfa,\bfb,\bfc)$
could be computed.
Specifically, we computed $q(\bfzero,\bfzero,\bfc)$ for \textit{all}
samples of a given size $c$ using the method described in \citet
{jenson2009G}. By weighting\vspace*{1pt} samples according to their sampling
probability, we may compute the distributions of $\ePS{M}$ and $\ePade
{M}$, and similarly the distributions of $\eaPS{M}$ and $\eaPade{M}$.
Table \ref{tabmain} summarizes the cumulative distributions of $\eaPS{M}$
and $\eaPade{M}$ for $\rho= 50$,
across all samples of size $c=20$ that are dimorphic at both loci. The
corresponding table for $\ePS{M}$ and $\ePade{M}$ was essentially
identical (not shown), with agreement usually to two decimal places.
This confirms that utilizing Approximation \ref
{approxrecursionapproximation} is
justified.

%
%
\begin{table}
\caption{Cumulative distribution\vspace*{-1pt}
$\Phi(x) = \bbP(\protect\eApprox^{(M)} <
x\%)$
(where $\protect\eApprox^{(M)}$ denotes either $\protect\eaPS{M}$ or
$\protect\eaPade{M}$) of the unsigned relative error of the partial sum
$\qPS{M}$ and
the corresponding Pad\'e approximants, for all samples of size $20$
dimorphic at both loci. Here, $\rho= 50$, and Approximation
\protect\ref{approxrecursionapproximation} is used}
\label{tabmain}
\begin{tabular*}{\tablewidth}{@{\extracolsep{\fill}}lcd{1.3}d{1.3}cccc@{}}
\hline
$\bolds{M}$ & \textbf{Type of sum} & \multicolumn{1}{c}{$\bolds{\Phi(1)}$} &
\multicolumn{1}{c}{$\bolds{\Phi(5)}$} & $\bolds{\Phi(10)}$ & $\bolds{\Phi(25)}$ &
$\bolds{\Phi(50)}$ & $\bolds{\Phi(100)}$ \\
\hline
\hphantom{0}0 & PS & 0.49\tabnoteref{ta} & 0.56\tabnoteref{ta} & 0.63 & 0.81 &
0.98 & 1.00 \\
& Pad\'e & 0.49 & 0.56 & 0.63 & 0.81 & 0.98 & 1.00 \\[4pt]
\hphantom{0}1 & PS & 0.51 & 0.74 & 0.87 & 0.99 & 1.00 & 1.00 \\
& Pad\'e & 0.59 & 0.77 & 0.84 & 0.91 & 0.98 & 0.99 \\[4pt]
\hphantom{0}2 & PS & 0.59 & 0.87 & 0.92 & 1.00 & 1.00 & 1.00 \\
& Pad\'e & 0.77 & 0.97 & 0.98 & 0.99 & 1.00 & 1.00 \\[4pt]
\hphantom{0}3 & PS & 0.57 & 0.86 & 0.97 & 1.00 & 1.00 & 1.00 \\
& Pad\'e & 0.91 & 0.96 & 0.98 & 0.99 & 1.00 & 1.00 \\[4pt]
\hphantom{0}4 & PS & 0.45 & 0.62 & 0.87 & 0.98 & 0.98 & 1.00 \\
& Pad\'e & 0.95 & 0.99 & 1.00 & 1.00 & 1.00 & 1.00 \\[4pt]
\hphantom{0}5 & PS & 0.30 & 0.50 & 0.61 & 0.76 & 0.79 & 0.97 \\
& Pad\'e & 0.98 & 1.00 & 1.00 & 1.00 & 1.00 & 1.00 \\[4pt]
10 & PS & 0.00 & 0.02 & 0.03 & 0.07 & 0.08 & 0.11 \\
& Pad\'e & 1.00 & 1.00 & 1.00 & 1.00 & 1.00 & 1.00 \\ [4pt]
& OTR & 0.25 & 0.36 & 0.65 & 0.89 & 0.90 & 1.00 \\
\hline
\end{tabular*}
\tabnotetext[\mbox{$\dagger$}]{ta}{These two values were misquoted in the text
of Jenkins and Song [(\protect\citeyear{jenson2009G}), page 1093]; this
table corrects them.}
\end{table}
%

Table \ref{tabmain} illustrates that the observations of Section \ref
{secexample} hold much more generally, as described below.
\begin{longlist}[(3)]
\item[(1)]
The error $\eaPS{M}$ for partial sums is not a monotonically
decreasing function of $M$; that is, the accuracy of $\qaPS{M}$
improves as one adds more terms up to a certain point, before quickly
becoming very inaccurate.
\item[(2)]
Empirically, the actual optimal truncation point for the parameter
settings we considered is at $M'=2$ or $M'=3$, which perform
comparably. Moreover, both provide consistently higher accuracy than
employing the OTR, which is a serious issue when we wish to use this
rule without external information about which truncation point really
is the most accurate.
\item[(3)]
Overall, using Pad\'e approximants is much more reliable. Note that the
accuracy of Pad\'e approximants continues to improve as we incorporate
more terms.
For a sample drawn at random, the probability that its Pad\'e
approximant is within $1\%$ of the true sampling distribution is
$1.00$, compared to $0.59$ or $0.57$ for truncating the partial sums,
respectively, at $M'=2$ or $M'=3$, and only $0.25$ for using the OTR.
\end{longlist}

For the remainder of this section, we focus on the accuracy of the
staircase Pad\'e approximants.
It was shown in \citet{jenson2009G} that the accuracy of the
partial sum $\qPS{2}$ increases with $\rho$, but (perhaps
surprisingly) decreases with increasing sample size. In Tables \ref
{tabrplot} and \ref{tabcplot}, we address the same issue for the
%
%
\begin{table}
\caption{Effect of $\rho$\vspace*{1pt} on the cumulative distribution
$\Phi(x) =
\bbP(\protect\eApprox_{\Pade}^{(M)} < x\%)$ of the unsigned relative
error of the Pad\'e approximants, for all samples of size $20$
dimorphic at both loci}
\label{tabrplot}
\begin{tabular*}{\tablewidth}{@{\extracolsep{\fill}}rd{3.0}cccccc@{}}
\hline
$\bolds{M}$ & \multicolumn{1}{c}{$\bolds\rho$} & $\bolds{\Phi(1)}$ & $\bolds{\Phi(5)}$ & $\bolds{\Phi(10)}$ &
$\bolds{\Phi(25)}$ & $\bolds{\Phi(50)}$ & $\bolds{\Phi(100)}$ \\
\hline
0 & 25 & 0.39 & 0.52 & 0.58 & 0.69 & 0.84 & 1.00 \\
& 50 & 0.49 & 0.56 & 0.63 & 0.81 & 0.98 & 1.00 \\
& 100 & 0.50 & 0.61 & 0.72 & 0.97 & 1.00 & 1.00 \\
& 200 & 0.54 & 0.72 & 0.95 & 0.99 & 1.00 & 1.00 \\[4pt]
1 & 25 & 0.51 & 0.62 & 0.75 & 0.85 & 0.94 & 0.96 \\
& 50 & 0.59 & 0.77 & 0.84 & 0.91 & 0.98 & 0.99 \\
& 100 & 0.74 & 0.91 & 0.95 & 0.98 & 0.99 & 1.00 \\
& 200 & 0.90 & 0.98 & 0.99 & 1.00 & 1.00 & 1.00 \\[4pt]
2 & 25 & 0.59 & 0.82 & 0.91 & 0.94 & 0.95 & 0.97 \\
& 50 & 0.77 & 0.97 & 0.98 & 0.99 & 1.00 & 1.00 \\
& 100 & 0.95 & 1.00 & 1.00 & 1.00 & 1.00 & 1.00 \\
& 200 & 1.00 & 1.00 & 1.00 & 1.00 & 1.00 & 1.00 \\[4pt]
3 & 25 & 0.64 & 0.91 & 0.95 & 0.96 & 0.98 & 1.00 \\
& 50 & 0.91 & 0.96 & 0.98 & 0.99 & 1.00 & 1.00 \\
& 100 & 0.99 & 1.00 & 1.00 & 1.00 & 1.00 & 1.00 \\
& 200 & 1.00 & 1.00 & 1.00 & 1.00 & 1.00 & 1.00 \\[4pt]
4 & 25 & 0.83 & 0.96 & 0.99 & 0.99 & 1.00 & 1.00 \\
& 50 & 0.95 & 0.99 & 1.00 & 1.00 & 1.00 & 1.00 \\
& 100 & 1.00 & 1.00 & 1.00 & 1.00 & 1.00 & 1.00 \\
& 200 & 1.00 & 1.00 & 1.00 & 1.00 & 1.00 & 1.00 \\[4pt]
5 & 25 & 0.82 & 0.94 & 0.98 & 0.99 & 1.00 & 1.00 \\
& 50 & 0.98 & 1.00 & 1.00 & 1.00 & 1.00 & 1.00 \\
& 100 & 1.00 & 1.00 & 1.00 & 1.00 & 1.00 & 1.00 \\
& 200 & 1.00 & 1.00 & 1.00 & 1.00 & 1.00 & 1.00 \\[4pt]
10 & 25 & 0.97 & 0.99 & 0.99 & 1.00 & 1.00 & 1.00 \\
& 50 & 1.00 & 1.00 & 1.00 & 1.00 & 1.00 & 1.00 \\
& 100 & 1.00 & 1.00 & 1.00 & 1.00 & 1.00 & 1.00 \\
& 200 & 1.00 & 1.00 & 1.00 & 1.00 & 1.00 & 1.00 \\
\hline
\end{tabular*}
\end{table}
Pad\'e approximant. Table \ref{tabrplot} confirms that accuracy increases
with $\rho$, as one might expect. Furthermore, it is also the case
that substantial accuracy is achievable even for moderate values of
$\rho$ (say, $\rho= 25$), provided that sufficiently many terms are
utilized in the construction of the Pad\'e approximant. For example,
when $M=5$ the probability that the Pad\'e approximant of a sample
drawn at random is within $5\%$ of the truth is $0.94$. Also very
encouraging is the pattern shown in Table \ref{tabcplot}. Provided that
%
%
\begin{table}
\caption{Effect of sample size $c$ on the cumulative distribution
$\Phi(x) = \bbP(\protect\eApprox_{\Pade}^{(M)} < x\%)$ of the
unsigned relative error of the Pad\'e approximants, for all samples of
size $c$ dimorphic at~both loci. Here~$\rho= 50$}
\label{tabcplot}
\begin{tabular*}{\tablewidth}{@{\extracolsep{\fill}}rccccccc@{}}
\hline
$\bolds{M}$ & $\bolds{c}$ & $\bolds{\Phi(1)}$ & $\bolds{\Phi(5)}$ & $\bolds{\Phi(10)}$ &
$\bolds{\Phi(25)}$ & $\bolds{\Phi(50)}$ & $\bolds{\Phi(100)}$ \\
\hline
0 & 10 & 0.58 & 0.67 & 0.72 & 0.96 & 1.00 & 1.00 \\
& 20 & 0.49 & 0.56 & 0.63 & 0.81 & 0.98 & 1.00 \\
& 30 & 0.44 & 0.50 & 0.58 & 0.72 & 0.91 & 1.00 \\[4pt]
1 & 10 & 0.72 & 0.90 & 0.94 & 0.97 & 1.00 & 1.00 \\
& 20 & 0.59 & 0.77 & 0.84 & 0.91 & 0.98 & 0.99 \\
& 30 & 0.53 & 0.69 & 0.76 & 0.85 & 0.95 & 0.97 \\[4pt]
2 & 10 & 0.94 & 1.00 & 1.00 & 1.00 & 1.00 & 1.00 \\
& 20 & 0.77 & 0.97 & 0.98 & 0.99 & 1.00 & 1.00 \\
& 30 & 0.62 & 0.89 & 0.95 & 0.98 & 0.99 & 0.99 \\[4pt]
3 & 10 & 0.98 & 1.00 & 1.00 & 1.00 & 1.00 & 1.00 \\
& 20 & 0.91 & 0.96 & 0.98 & 0.99 & 1.00 & 1.00 \\
& 30 & 0.81 & 0.92 & 0.94 & 0.98 & 0.99 & 1.00 \\[4pt]
4 & 10 & 1.00 & 1.00 & 1.00 & 1.00 & 1.00 & 1.00 \\
& 20 & 0.95 & 0.99 & 1.00 & 1.00 & 1.00 & 1.00 \\
& 30 & 0.90 & 0.99 & 0.99 & 1.00 & 1.00 & 1.00 \\[4pt]
5 & 10 & 1.00 & 1.00 & 1.00 & 1.00 & 1.00 & 1.00 \\
& 20 & 0.98 & 1.00 & 1.00 & 1.00 & 1.00 & 1.00 \\
& 30 & 0.86 & 0.99 & 0.99 & 1.00 & 1.00 & 1.00 \\[4pt]
10 & 10 & 1.00 & 1.00 & 1.00 & 1.00 & 1.00 & 1.00 \\
& 20 & 1.00 & 1.00 & 1.00 & 1.00 & 1.00 & 1.00 \\
& 30 & 0.99 & 1.00 & 1.00 & 1.00 & 1.00 & 1.00 \\
\hline
\end{tabular*}
\end{table}
sufficiently many terms are used, the accuracy of the Pad\'e
approximant is only slightly affected by looking at larger sample
sizes. For example, when $M=5$ the probability that the Pad\'e
approximant of a randomly drawn sample of size $30$ is within $5\%$ of
the truth is $0.99$, compared with $1.00$ for a sample of size $20$.
This loss in accuracy is much less severe than the corresponding loss
in accuracy of the partial sums;
for $\qPS{M}$, the highest accuracy is achieved for $M = 2$, in which
case which the corresponding loss in accuracy is from $0.87$ when $c =
20$, to $0.70$ when $c = 30$.

\section{Proofs}
\label{secproofs}
In this section, we provide proofs of the results presented earlier.

\subsection{\texorpdfstring{Proof of Theorem \protect\ref{thmg}}{Proof of Theorem 3.1}}
\label{secproofg}
For an infinitesimal generator $\sA$ of a diffusion process on the
state space $\Omega$ and a twice continuously differentiable function
$h\dvtx\Omega\to\bbR$ with compact support, it is well known that
\[
\bbE[\sA h(\bfX)] = 0,
\]
where expectation is with respect to the stationary distribution of
$\bfX$. Apply this result to the generator $\widetilde{\sL}$ shown
in (\ref{eqgenneutral2}) and monomial $G^{(m)}(\widetilde{{\bfmath
{x}}};\bfa,\bfb,\bfr)$ shown in~(\ref{eqgm}). This provides a
linear equation relating the expectations $\bbE[G^{(m+1)}(\tbfX;\bfa
'$, $\bfb',\bfr')]$, $\bbE[G^{(m)}(\tbfX;\bfa',\bfb',\bfr')]$,
$\bbE[G^{(m-1)}(\tbfX;\bfa',\bfb',\bfr')]$ and\break $\bbE
[G^{(m-2)}(\tbfX;\bfa',\bfb'$, $\bfr')]$, each appearing with various
different arguments $(\bfa',\bfb',\bfr')$ depending on $(\bfa,\bfb
,\bfr)$; we omit the simple but algebraically lengthy details. Now,
for these four choices of $m$ we substitute the proposed expansion
(\ref{eqexpansion}). If $m>0$, then compare the coefficients of $\rho
^{1-{u/2}}$ in the resulting expression; if $m=0$, then
compare the coefficients of $\rho^{-{ u/2}}$. We then obtain
the following:\vspace*{8pt}

\mbox{}\hphantom{i}(i) If $m=0$ and $ u=0$, then the resulting expression is
%
%
\begin{eqnarray}
\label{eqcomputeg0}
&&
[a(a-1+\tA) + b(b-1+\tB)]g_0^{(0)}(\bfa,\bfb,\bfzero)
\nonumber\\[-2pt]
&&\qquad= \sum_i a_i(a_i - 1) g_0^{(0)}(\bfa-\bfe_i,\bfb,\bfzero) + \sum
_j b_j(b_j - 1) g_0^{(0)}(\bfa,\bfb-\bfe_j,\bfzero)
\nonumber\\[-9pt]\\[-9pt]
&&\qquad\quad{} + \tA\sum_{i,k} a_i P_{ki}^\A g_0^{(0)}(\bfa-\bfe_i + \bfe
_k,\bfb,\bfzero) \nonumber\\[-2pt]
&&\qquad\quad{} + \tB\sum_{j,l} b_j P_{lj}^\B g_0^{(0)}(\bfa,\bfb-\bfe_j +
\bfe_l,\bfzero)\nonumber
\end{eqnarray}
with boundary conditions given by (\ref{eqnewboundaries}). This is the
sum of two copies of a~familiar recursion [\citet{gritav1994TPB}]
for the sampling distribution of a single locus, one for locus $A$ and
one for locus $B$. In our notation the solution is, therefore,
\[
g^{(0)}_0(\bfa,\bfb,\bfzero) = q^\A(\bfa)q^\B(\bfb).
\]
%

(ii) If $m>0$, then the resulting expression is
%
%
\begin{eqnarray}
\label{eqcomputeg}
\hspace*{-4pt}&&mg_ u^{(m)}(\bfa,\bfb,\bfr)  \nonumber
\\[-2pt]
\hspace*{-4pt}&&\qquad=\sum_{i,j}\biggl\{ r_{ij}(r_{ij}-1)g_ u
^{(m-2)}(\bfa+\bfe_i,\bfb+\bfe_j,\bfr-2\bfe_{ij})
\nonumber\\[-2pt]
\hspace*{-4pt}&&\qquad\hphantom{=\sum_{i,j}\biggl\{}{} -\sum_l r_{ij}(r_{il}-\delta_{jl})g_ u^{(m-2)}(\bfa+\bfe
_i,\bfb+\bfe_j+\bfe_l,\bfr-\bfe_{ij}-\bfe_{il})\nonumber\\[-2pt]
\hspace*{-4pt}&&\qquad\hphantom{=\sum_{i,j}\biggl\{}{} - \sum_k r_{ij}(r_{kj} - \delta_{ik})g_ u^{(m-2)}(\bfa+\bfe
_i+\bfe_k,\bfb+\bfe_j,\bfr-\bfe_{ij}-\bfe_{kj})\nonumber\\[-2pt]
\hspace*{-4pt}&&\qquad\hphantom{=\sum_{i,j}\biggl\{}{} + \sum_{k,l} r_{ij}(r_{kl}-\delta_{ik}\delta_{jl})g_ u
^{(m-2)}(\bfa+\bfe_i+\bfe_k,\bfb+\bfe_j+\bfe_l,\bfr-\bfe
_{ij}-\bfe_{kl})\nonumber\\[-2pt]
\hspace*{-4pt}&&\qquad\hphantom{=\sum_{i,j}\biggl\{}{} + r_{ij}(r_{ij}-1)g_{ u-1}^{(m-1)}(\bfa,\bfb,\bfr-\bfe
_{ij})\nonumber\\[-2pt]
\hspace*{-4pt}&&\qquad\hphantom{=\sum_{i,j}\biggl\{}{} -2 r_{ij}(r_{i\cdot}-1)g_{ u-1}^{(m-1)}(\bfa,\bfb+\bfe_j,\bfr
-\bfe_{ij})\nonumber\\[-2pt]
\hspace*{-4pt}&&\qquad\hphantom{=\sum_{i,j}\biggl\{}{} -2 r_{ij}(r_{\cdot j}-1)g_{ u-1}^{(m-1)}(\bfa+\bfe_i,\bfb,\bfr
-\bfe_{ij})\nonumber\\[-2pt]
\hspace*{-4pt}&&\qquad\hphantom{=\sum_{i,j}\biggl\{}{} +2\sum_{k,l}r_{kj}(r_{il}-\delta_{ik}\delta_{jl})g_{ u
-1}^{(m-1)}(\bfa+\bfe_k,\bfb+\bfe_l,\bfr-\bfe_{kj}-\bfe
_{il}+\bfe_{ij})\nonumber\\[-2pt]
\hspace*{-4pt}&&\qquad\hspace*{101pt}\hphantom{=\sum_{i,j}\biggl\{}{} +2(m-1)r_{ij}g_{ u-1}^{(m-1)}(\bfa+\bfe_i,\bfb+\bfe
_j,\bfr-\bfe_{ij}) \biggr\}\nonumber\\[-2pt]
\hspace*{-4pt}&&\qquad\quad{} +\sum_i a_i(a_i+2r_{i\cdot}-1)g_{ u-2}^{(m)}(\bfa-\bfe_i,\bfb
,\bfr)\\[-2pt]
\hspace*{-4pt}&&\qquad\quad{} +\sum_j b_j(b_j+2r_{\cdot j}-1) g_{ u-2}^{(m)}(\bfa,\bfb-\bfe
_j,\bfr) \nonumber\\[-2pt]
\hspace*{-4pt}&&\qquad\quad{} -2\sum_{i,j}\sum_k a_i r_{kj} g_{ u-2}^{(m)}(\bfa-\bfe_i+\bfe
_k,\bfb,\bfr-\bfe_{kj}+\bfe_{ij})\nonumber\\[-2pt]
\hspace*{-4pt}&&\qquad\quad{} -2\sum_{i,j}\sum_l b_j r_{il} g_{ u-2}^{(m)}(\bfa,\bfb-\bfe
_j+\bfe_l,\bfr-\bfe_{il}+\bfe_{ij})\nonumber\\[-2pt]
\hspace*{-4pt}&&\qquad\quad{} +\tA\sum_{i,k}P_{ki}^\A\biggl[a_i g_{ u-2}^{(m)}(\bfa-\bfe
_i+\bfe_k,\bfb,\bfr) \nonumber\\[-2pt]
\hspace*{-4pt}&&\qquad\quad\hphantom{{}+\tA\sum_{i,k}P_{ki}^\A\biggl[}{} +\sum_j r_{ij}g_{ u
-2}^{(m)}(\bfa,\bfb,\bfr-\bfe_{ij}+\bfe_{kj})\biggr]\nonumber\\[-2pt]
\hspace*{-4pt}&&\qquad\quad{} +\tB\sum_{j,l}P_{lj}^\B\biggl[b_j g_{ u-2}^{(m)}(\bfa,\bfb
-\bfe_j+\bfe_l,\bfr) \nonumber\\[-2pt]
\hspace*{-4pt}&&\qquad\quad\hphantom{{}+\tB\sum_{j,l}P_{lj}^\B\biggl[}{} + \sum_i
r_{ij}g_{ u-2}^{(m)}(\bfa,\bfb,\bfr-\bfe_{ij}+\bfe_{il})
\biggr]\nonumber\\[-2pt]
\hspace*{-4pt}&&\qquad\quad{} - [(a+m)(a+m+\tA-1)\nonumber\\[-2pt]
\hspace*{-4pt}&&\qquad\quad\hphantom{{}- [}{} + (b+m)(b+m+\tB-1)  - m(m-3)]g_{ u-2}^{(m)}(\bfa,\bfb,\bfr)\nonumber\\[-2pt]
\hspace*{-4pt}&&\qquad\quad{} + 2\sum_{i,j}a_i b_j g_{ u-3}^{(m+1)}(\bfa-\bfe_i,\bfb-\bfe
_j,\bfr+\bfe_{ij}).\nonumber
\end{eqnarray}
Equation (\ref{eqcomputeg}) relates $g_ u^{(m)}(\bfa,\bfb,\bfr)$
to the known expressions $g_ u^{(m-2)}$, $g_{ u-1}^{(m-1)}$, $g_{ u
-2}^{(m)}$ and $g_{ u-3}^{(m+1)}$, as claimed.\vspace*{2pt}

(iii) If $m=0$ and $ u\geq1$, then the resulting expression is
%
%
\begin{eqnarray}
\label{eqcomputeg2}
&&
[a(a+\theta_\A-1) + b(b+\theta_\A-1)]g^{(0)}_ u(\bfa
,\bfb,\bfzero) \nonumber\\
&&\qquad = \sum_i a_i(a_i-1)g^{(0)}_ u(\bfa-\bfe_i,\bfb,\bfzero) + \sum
_j b_j(b_j-1) g^{(0)}_ u(\bfa,\bfb-\bfe_j,\bfzero)\nonumber\\
&&\qquad\quad{} + \tA\sum_{i,k} a_i P^\A_{ki} g^{(0)}_ u(\bfa-\bfe_i+\bfe
_k,\bfb,\bfzero) \\
&&\qquad\quad{} + \tB\sum_{j,l} b_j P^\B_{lj} g^{(0)}_ u(\bfa,\bfb-\bfe
_j+\bfe_l,\bfzero) \nonumber\\
&&\qquad\quad{} + 2\sum_{i,j}a_i b_j g^{(1)}_{ u-1}(\bfa-\bfe_i,\bfb-\bfe
_j,\bfe_{ij})\nonumber
\end{eqnarray}
with boundary conditions (\ref{eqnewboundaries}). Hence, this provides
a recursion relation for $g^{(0)}_u(\bfa,\bfb,\bfzero)$ when
$g^{(1)}_{u-1}$ is known.

\subsection{\texorpdfstring{Proof of Corollary \protect\ref{corodd}}{Proof of Corollary 3.1}}
\label{secproofodd}
For the base case $\ell=1$, note that if $m=1$ and $ u=0$ then (\ref
{eqcomputeg}) simplifies to
\[
g^{(1)}_0(\bfa,\bfb,\bfr) = 0,
\]
and hence, if $m=0$ and $ u=1$ then (\ref{eqcomputeg2}) simplifies to
%
%
\begin{eqnarray}
\label{eqnullrecursion}
&&[a(a+\theta_\A-1) + b(b+\theta_\A-1)]g^{(0)}_1(\bfa,\bfb
,\bfzero) \nonumber\\
&&\qquad= \sum_i a_i(a_i-1)g^{(0)}_1(\bfa-\bfe_i,\bfb,\bfzero) + \sum_j
b_j(b_j-1) g^{(0)}_1(\bfa,\bfb-\bfe_j,\bfzero)\nonumber\\[-8pt]\\[-8pt]
&&\qquad\quad{} + \tA\sum_{i,k} a_i P^\A_{ki} g^{(0)}_1(\bfa-\bfe_i+\bfe
_k,\bfb,\bfzero) \nonumber\\
&&\qquad\quad{} + \tB\sum_{j,l} b_j P^\B_{lj} g^{(0)}_1(\bfa,\bfb-\bfe_j+\bfe
_l,\bfzero)\nonumber
\end{eqnarray}
with boundary conditions (\ref{eqnewboundaries}). This has solution
\[
g^{(0)}_1(\bfa,\bfb,\bfzero) = 0,
\]
which is unique [since $q(\bfa,\bfb,\bfc)$ is unique]. This
completes $\ell=1$. Now suppose inductively that $g_ u^{(m)}(\bfa
,\bfb,\bfr) = 0$ for every $m$, $ u$ such that $m+ u=\ell-2$ where
$\ell$ is a fixed odd number greater than $1$. Then for $m$, $ u$
such that $m+ u=\ell$, (\ref{eqcomputeg})~becomes
\[
g_ u^{(m)}(\bfa,\bfb,\bfr) = 0
\]
as required.

\subsection{\texorpdfstring{Proof of Lemma \protect\ref{lemrational}}{Proof of Lemma 4.1}}
\label{secproofrational}
In what follows, define the {\sl length} of a sample configuration
$(\bfa,\bfb,\bfc)$ to be $a+b+2c$.
Under a neutral, finite-alleles model, the probability of
a sample with length $\delta$ satisfies a closed system of
equations\vadjust{\goodbreak}
[e.g., see equation (5) of \citet{jenson2009G}] which can be
expressed in matrix form:
\[
\bfM\bfq= \bfv,
\]
where $\bfq$ is a vector composed of the probabilities of samples of
length less than or equal to $\delta$, $\bfv$ is a constant vector of
the same dimension as $\bfq$ and~$\bfM$ is an invertible matrix
(since the solution to this equation is unique). The entries of $\bfM$
and $\bfv$ are rational functions of $\rho$, and hence, $\bfq= \bfM
^{-1}\bfv$ is a~vector each of whose entries is a~rational function of
$\rho$.

Let $U_0$ denote the degree of the numerator, and $V_0$ the degree of
the denominator. If $U_0 > V_0$, then $q(\bfa,\bfb,\bfc)$ becomes
unbounded as $\rho\to\infty$, while if $V_0 > U_0$ then $q(\bfa
,\bfb,\bfc) \to0$ as $\rho\to\infty$. But we know that $q(\bfa
,\bfb,\bfc)$ is a~probability, and hence, bounded. Moreover, it has
support over all samples of a~fixed size, since we assume that $\bfP
^\A$ and $\bfP^\B$ are irreducible. Thus, to ensure $\lim_{\rho\to
\infty} q(\bfa,\bfb,\bfc) \in(0,1)$ we must have $U_0 = V_0$. By a
similar\vspace*{1pt} argument as $\rho\to0$, we must have that the coefficients of
$\rho^0$ are nonzero both in the numerator and in the denominator. We
can, therefore, divide the numerator and denominator by $\rho^{U_0}$
to obtain a rational function of $1/\rho$ whose degree in the
numerator and denominator are both $U_0$ (\mbox{$=$}$V_0$).

\subsection{\texorpdfstring{Proof of Theorem \protect\ref{thmconvergence}}{Proof of Theorem 4.1}}
\label{secproofconvergence}
This is an application of Theorem 1.4.4 of \citet{bakgra1996},
which we spell out for completeness. By Lem\-ma~\ref{lemrational},
$q(\bfa,\bfb,\bfc)$ is a rational function of $1/\rho$ and is
analytic at $\rho= \infty$ with Taylor series~(\ref{eqmainexpansion}).
Denote the degree of its numerator and denominator
by $U_0$. Then, $q(\bfa,\bfb,\bfc)$ has $U_0 + U_0 + 1$ independent
coefficients determined by the first $U_0 + U_0 + 1$ terms of its
Taylor series expansion. Thus, provided $U \geq U_0$ and $V \geq U_0$,
by the definition of the $[U/V]$ Pad\'e approximant, it must coincide
uniquely with $q(\bfa,\bfb,\bfc)$.

\subsection{\texorpdfstring{Proof of Theorem \protect\ref{thmmain-s}}{Proof of Theorem 5.1}}
\label{secproofmain-s}
This is simply an application of Theorem \ref{thmg} applied to the
generator for the diffusion under selection, $\widetilde{\sL}_\ts$,
rather than~$\widetilde{\sL}$, and so we just summarize the procedure.

The change in generator results in slight modifications to the
relationships between the $g^{(m)}_ u(\bfa,\bfb,\bfr)$ in order to
obtain the relationships between the $h^{(m)}_ u(\bfa,\bfb,\bfr)$
for each $m$ and $ u$:
\begin{longlist}[(1)]
\item[(1)] $h^{(0)}_0(\bfa,\bfb,\bfzero)$ satisfies (\ref{eqcomputeg0})
(replacing each $g^{(0)}_0$ with $h^{(0)}_0$), but with extra terms
\[
{}+ \sum_{i,k} a_i\biggl[\gs_{ik} h_0^{(0)}(\bfa+\bfe_k,\bfb
,\bfzero) - \sum_{k'}\gs^\A_{kk'}h_0^{(0)}(\bfa+\bfe_k+\bfe
_{k'},\bfb,\bfzero)\biggr]
\]
on the right-hand side. The solution is
\[
h_0^{(0)}(\bfa,\bfb,\bfzero) = q^\A_\ts(\bfa)q^\B(\bfb).
\]

\item[(2)] For $m>0$, $h^{(0)}_0(\bfa,\bfb,\bfzero)$ satisfies (\ref
{eqcomputeg}) but with extra terms
%
%
\begin{eqnarray}
\label{eqcomputeg-s}
&&{}+ \sum_{i,k}(a_i + r_{i\cdot})\gs_{ik}^\A h^{(m)}_{ u-2}(\bfa
+\bfe_k,\bfb,\bfr) \nonumber\\
&&\qquad{}- (a+m)\sum_{k,k'}\gs_{kk'}^\A h^{(m)}_{ u-2}(\bfa+\bfe_k+\bfe
_{k'},\bfb,\bfr) \\
&&\qquad{}- \sum_{i,j,k,k'} r_{ij} \gs_{kk'}^\A h^{(m)}_{ u-2}(\bfa+\bfe
_i+\bfe_k,\bfb,\bfr-\bfe_{ij}+\bfe_{k'j})
\nonumber
\end{eqnarray}
on the right-hand side.
\item[(3)] For $m=0$ and $ u\geq1$, $h^{(0)}_0(\bfa,\bfb,\bfzero)$
satisfies (\ref{eqcomputeg}) but with extra terms
\begin{eqnarray*}
&&{}+ \sum_{i,k}\bigl[a_i\gs_{ik}^\A h^{(0)}_ u(\bfa+\bfe_k,\bfb
,\bfzero) - a\gs_{ik}^\A h^{(0)}_ u(\bfa+\bfe_i+\bfe_k,\bfb
,\bfzero)\bigr] \\
&&\qquad{}+ \sum_{i,j,k} b_j\gs_{ik}^\A h^{(1)}_{ u-1}(\bfa+\bfe_k,\bfb
-\bfe_j,\bfe_{ij})
\end{eqnarray*}
on the right-hand side.
\end{longlist}
Using these equations to evaluate $h^{(m)}_ u(\bfa,\bfb,\bfr)$ on
levels $\ell=0,1,\ldots,4$ provides expressions for the nonneutral
versions of $q_0(\bfa,\bfb,\bfc)$, $q_1(\bfa,\bfb,\bfc)$ and
$q_2(\bfa,\bfb,\bfc)$. Those for $q_0(\bfa,\bfb,\bfc)$ and
$q_1(\bfa,\bfb,\bfc)$ in terms of the relevant one-locus sampling
distributions are unchanged, while the new generator makes some minor
modifications to the expression for $q_2(\bfa,\bfb,\bfc)$. The
analytic part of this term, $\sigma_\ts(\bfa,\bfb,\bfc)$, is
easily calculated from $h^{(4)}_0$, $h^{(3)}_1$, $h^{(2)}_2$ and $h^{(1)}_3$,
while the recursive part, $q_{2,\ts}(\bfa+\bfc_\A,\bfb+\bfc_\B
,\bfzero)$, follows from $h^{(0)}_4$.


\section*{Acknowledgments}
We gratefully acknowledge Anand Bhaskar for implementing the algorithm
stated in Theorem \ref{thmg} and for checking some of our formulas. We
also thank Matthias Steinr\"ucken and Paul Fearnhead for useful
discussion.


%

%
\printaddresses

\end{document}